\documentclass[a4paper, 12pt]{article}
\usepackage{enumerate, theorem}
\usepackage{amsmath, amsfonts, amssymb}
\usepackage[height=22.5cm, width=15cm]{geometry}
\usepackage[all]{xy}
\usepackage{mathrsfs, yfonts}

\DeclareMathOperator{\Sym}{Sym}

\DeclareMathOperator{\C}{\mathbb{C}}

\newcommand{\parag}[1]{\paragraph{\sc{#1.}} }

\newtheorem{thm}{Theorem}[subsection]
\newtheorem{defn}[thm]{Definition}
\newtheorem{cor}[thm]{Corollary}
\newtheorem{prop}[thm]{Proposition}
\newtheorem{lemma}[thm]{Lemma}

\begin{document}

\title{The sheaf $\alpha_{X}^{\bullet}$ }
\date{20/07/17}

 \author{Daniel Barlet\footnote{Institut Elie Cartan, G\'eom\`{e}trie,\newline
Universit\'e de Lorraine, CNRS UMR 7502   and  Institut Universitaire de France.}.}

\maketitle

\tableofcontents

\parag{Abstract} We introduce in a reduced complex space, a ``new coherent sub-sheaf'' of the sheaf $\omega_{X}^{\bullet}$ which has the ``universal pull-back property'' for any holomorphic map, and which is in general bigger than the usual sheaf of holomorphic differential forms $\Omega_{X}^{\bullet}/torsion$. We show that the meromorphic differential forms which are sections of this sheaf satisfy integral dependence equations over the symmetric algebra of the sheaf $\Omega_{X}^{\bullet}/torsion$. This sheaf  $\alpha_{X}^{\bullet}$ is also closely related to the normalized Nash transform.\\
We also show that these $q-$meromorphic differential forms are locally  square-integrable on {\bf any} $q-$dimensional cycle in $X$ and that the corresponding functions obtained by integration on an analytic family of $q-$cycles are locally bounded and  locally continuous on the complement of closed analytic subset.

\parag{AMS classification} 32 C 15 - 32 C 30 - 32 S xx - 32 S 45.

\parag{Key words} Meromorphic differential forms on a singular space- Universal pull-back property- Normalized Nash transform- Integral dependence equation for differential forms.

\section*{Introduction}

In this article we discuss the following question: given a reduced complex space $X$ the normalization of $X$ consists in building a proper modification $\nu : \tilde{X} \to X$ such that meromorphic locally bounded functions on $X$ becomes holomorphic after pull-back to $\tilde{X}$. Moreover this process gives a desingularization process for curves, that is to say for $X$ of  pure dimension 1.\\
It seems then natural to define an analogous process for meromorphic locally bounded differential forms. The main trouble is to define what means a ``locally bounded'' for a meromorphic differential form of positive  degree  on a reduced complex space. To define this notion is the purpose of this paper. Of course, this does not lead to a simple proof of a  desingularization process for a reduced complex space, but we will show that the natural  process associated to ``normalization of differential meromorphic forms'' is simply the classical {\bf normalized Nash transform}, and it is an old (an probably very difficult) conjecture that this process leads to a desingularization. We hope that the introduction of this ``new sheaf'' $\alpha_{X}^{\bullet}$ will be useful in that direction.\\
But in fact, the main reason to introduce this sheaf is the look for the {\bf ``universal pull-back property''} which means to define a coherent sheaf of meromorphic differential forms which admits a natural pull-back for {\bf any} holomorphic map between reduced complex spaces and which is ``maximal'' with this property. Note that if we only consider complex manifolds the sheaf $\Omega_{X}^{\bullet}$ has this property, but we will show that this is no longer maximal when $X$ admits singularities.\\
Our main result is the theorem \ref{pull-back} (and its precise formulation \ref{pull-back bis})  giving the ``universal pull-back property'' for these sheaves. We obtain also two other results which may be useful: 
\begin{itemize}
\item The fact that  for any section $\alpha$  of the sheaf $\alpha_{X}^{q}$ the form $\alpha\wedge \bar \alpha$ is locally integrable on {\bf any} holomorphic cycle of dimension $q$ and also the local boundness and the ``generic'' continuity of  such an integral  when the $q-$cycle moves in an analytic family (see theorem \ref{int. param.});
\item The existence of a local integral dependence equation for a section of $\alpha_{X}^{q}$ over the symmetric algebra of the sheaf $\Omega_{X}^{q}/torsion$ (see proposition \ref{int. dep.}). 
\end{itemize}
We conclude this article by computing some simple examples showing that the sheaf $\alpha_{X}^{\bullet}$ may be different from other classical sheaves of meromorphic differential forms which are used on singular complex spaces.

\section{Universal pull-back for $\Omega_{X}^{\bullet}/torsion$}

\parag{Conventions and Notations}\begin{itemize}
\item  Any complex space will be reduced in the sequel.
\item The graded sheaf $\Omega_{X}^{\bullet}$ on the complex space $X$ is locally defined in a local embedding of $X$ in an open set $ U \subset \C^{N}$ as the quotient
$$ \Omega_{\C^{N}}^{\bullet}\big/ \mathcal{I}_{X}.\Omega_{\C^{N}}^{\bullet} + d\mathcal{I}_{X}\wedge \Omega_{\C^{N}}^{\bullet-1}$$
where $\mathcal{I}_{X}$ is the reduced ideal of $X$ in $U$.
\item On a complex space $X$ we shall consider the graded sheaf $L_{X}^{\bullet}$ which is the direct image of the sheaf of holomorphic forms on a desingularization of $X$, and the graded sheaf \ $\omega_{X}^{\bullet} \subset j_{*}j^{*}\Omega_{X}^{\bullet}$, where $j : X \setminus S \to X$ is the inclusion of the regular part of $X$, which is the sheaf of $\bar\partial-$closed  $(\bullet, 0)-$currents on $X$ modulo these currents which have support in $S$ (see [B.78]). Then we have the natural inclusions of graded coherent sheaves
$$ \Omega_{X}^{\bullet}\big/torsion \subset L_{X}^{\bullet} \subset \omega_{X}^{\bullet} .$$
Note that in the case where the codimension of the singular set in $X$ is at least equal to $2$, we have the equality $\omega_{X}^{\bullet} = j_{*}j^{*}\Omega_{X}^{\bullet}$.
\item When $X$ is  not irreducible, let $X := \cup_{i \in I} X_{i}$ be the decomposition of $X$ in its irreducible components. In this case a desingularization of $X$ will be given by the disjoint union $\tau : \tilde{X} \to X$ of desingularizations $\tau_{i} : \tilde{X}_{i} \to X_{i} \hookrightarrow X$, where $\tilde{X}$ is the disjoint union of the connected  complex manifolds $\tilde{X}_{i}$.\\
It is easy to see that in this situation the sheaf $L_{X}^{\bullet}$ is the direct sum of the direct images in $X$ of the sheaves $L_{X_{i}}^{\bullet}$. This direct sum is locally finite on $X$ as the irreducible components of $X$ are locally finite on $X$.
\end{itemize}

It is clear that the sheaf  \ $\Omega_{X}^{\bullet}$ on a complex space has the ``universal pull-back'' property because the pull-back of a holomorphic form by a holomorphic map $f : Y \to X$  is well defined and functorial in the map $f$.\\

To begin we shall prove that the sheave $\Omega_{X}^{\bullet}/torsion$  still have this ``universal pull-back'' property.

\begin{prop}\label{torsion}
Let $X$ be a reduced complex space and consider a torsion holomorphic $p-$form $\alpha$ on $X$ ( meaning that $\alpha$ vanishes at smooth points in $X$). Let $Z$ be an analytic subset in $X$. Then the $p-$holomorphic form induced by $\alpha$ on $Z$ is again a torsion form on $Z$.
\end{prop}

\parag{Proof} Without any lost of generality we can assume that $Z$ is irreducible. Let $S$ the singular set of $X$. If $Z$ is not contained in $S$ the result is obvious. Also if the dimension of $Z$ is less that $p$ the conclusion is again obvious. So let $dim \, Z = p + q$ with $q \geq 0$ and let $Z'$ the dense open set of smooth points $x$ in $Z$   for which the multiplicity of $x$ in $X$ is minimal. It is enough to show that the restriction of $\alpha$ to $Z'$ vanishes. As the problem is local on $Z'$, we can assume that we have an open neighbourhood $X'$ of $x_{0}$ in $X$ and  a local parametrisation  $\pi : X' \to U$ on a polydisc $U$ of $\C^{n}$ with the following properties:
\begin{enumerate}[i)]
\item $\pi(x_{0}) = 0$.
\item $U = V \times W$ where  $V$ and $W$ are polydiscs  with center $0$ respectively in  $\C^{p+q}$  and  $\C^{n-p-q}$.
\item $Z'' := Z' \cap X'  = \pi^{-1}(V\times \{0\})$ set theoretically and $\pi : Z'' \to V\times \{0\}$ is an isomorphism.
\end{enumerate}
Define the analytic family  of $(p+q)-$cycles $(Z_{w})_{w \in W}$ in $X'$ parametrized by $W$ by letting $Z_{w} := \pi^{*}(V \times \{w\})$, where the pull-back by $\pi$ is taken in the sense of cycles\footnote{This means that if $f : U \to \Sym^{k}(X')$ is the holomorphic map classifying the fibers of $\pi$, the cycle $Z_{w}$ is the cycle-graph of the analytic  family of $k-$tuples in $X'$ defined by the restriction of $f$ to $V \times \{w\}$.}. Then, if $k$ is the degree of $\pi$ (which is the multiplicity in $X$ of each point in $Z''$) we have $Z_{0} = k.Z''$ as a cycle in $X'$. Remark that for $w$ generic in $W$ the intersection of the cycle $Z_{w}$ with the ramification set of $\pi$ has no interior point in $Z_{w}$ which is a reduced cycle. So the restriction of the holomorphic form $\alpha$ to $Z_{w}$ for $w$ generic is a torsion form.
\\
Now choose a non negative continuous function with compact support $\rho$ on $X'$ and a holomorphic $q-$form $\beta$ on $X'$ and define the function on $W$
$$ \varphi : W \to \mathbb{R}^{+}, \quad w \mapsto  \varphi(w) := \int_{Z_{w}} \ \rho.(\alpha\wedge \beta)\wedge \overline{(\alpha\wedge \beta)}.$$
It is a continuous function (see [B-M 1] ch.IV) and it vanishes for $w$ generic in  $W$ as $\alpha$ generically vanishes on $Z_{w}$ for such a $w$. Then it vanishes for $w = 0$ and this shows that the restriction of $\alpha$ to an open dense subset of $Z'$ vanishes.$\hfill \blacksquare$\\

\begin{cor}\label{torsion gen.}
Consider a holomorphic map $f : X \to Y$ where  $X$ and $Y$ are reduced complex spaces.Then, if $\alpha$ is a  $p-$holomorphic form on $Y$ which is a  torsion form, the $p-$holomorphic form  $f^{*}(\alpha)$ is a torsion form on $X$. 
\end{cor}

\parag{Proof} It is enough to consider the case where $X$ is a connected complex manifold. Let $X'$ be the open dense subset of $X$ where $f$ has maximal rank. On $X'$ the map $f$ is locally a submersion on a locally closed complex sub-manifold of $Y$ and the previous proposition applies to show that the pull-back of $\alpha$ on this locally closed sub-manifold vanishes. So the holomorphic 
 form $f^{*}(\alpha)$ vanishes on $X'$. Then it is a torsion form on $X$.$\hfill \blacksquare$
 
 \parag{Conclusion} \begin{itemize}
 \item The usual pull-back for holomorphic differential forms induced a natural pull-back for the sheaf $\Omega_{X}^{\bullet}/torsion$ by any holomorphic map between  reduced complex spaces.
 \end{itemize}
 
 \begin{defn}\label{double star} Let $f : X \to Y$ a holomorphic map between two reduced  complex spaces. We have a natural graded pull-back $\mathcal{O}_{X}-$morphism
 \begin{equation*} 
  f^{*} : f^{*}(\Omega_{Y}^{\bullet}/torsion) \to \Omega_{X}^{\bullet}/torsion \tag{*}
  \end{equation*}
  We shall denote $f^{**}(\Omega_{Y}^{\bullet})$ the image of this graded sheaf morphism.\\
  We shall also denote $f^{**}(\mathcal{G})$ for any sub-sheaf $\mathcal{G}$ of $\Omega_{Y}^{\bullet}/torsion$ its image by the morphism $f^{*}$ above (or also when $\mathcal{G}$ is a subs-sheaf of $\Omega^{\bullet}_{Y}$).
  \end{defn}
  So, by definition, $f^{**}(\Omega_{Y}^{\bullet})$ (and more generally $f^{**}(\mathcal{G})$) is a sub-sheaf of the sheaf $\Omega_{X}^{\bullet}/torsion$, so it has no $\mathcal{O}_{X}-$torsion.\\
  
  \begin{lemma}\label{pull-back 0}
  Let $f : X \to Y$ and $g : Y \to Z$ two holomorphic maps between reduced complex spaces. Then we have equality of the sub-sheaves $f^{**}(g^{**}(\mathcal{H}))$ and $(g\circ f)^{**}(\mathcal{H})$ for any sub-sheaf $\mathcal{H}$ of the sheaf $\Omega_{Z}^{\bullet}/torsion$.
  \end{lemma}
 
 \parag{Proof} The pull-back by $g$ gives a morphism
 $$ g^{*}(\Omega_{Z}^{\bullet}/torsion) \to \Omega_{Y}^{\bullet}/torsion $$
 with image $g^{**}(\Omega_{Z}^{\bullet}/torsion)$ and the pull-back by $f$ gives a morphism
 $$ f^{*}(g^{*}(\Omega_{Z}^{\bullet}/torsion)) \to f^{*}(\Omega_{Y}^{\bullet}/torsion) $$
 which, by right-exactness of the tensor product, is surjective on $f^{*}(g^{**}(\Omega_{Z}^{\bullet}/torsion))$. Then we have the following commutative diagram
 $$\xymatrix{f^{*}(g^{*}(\Omega_{Z}^{\bullet}/torsion)) \ar[r]^{\alpha} \ar[d]_{\simeq} & f^{*}(g^{**}(\Omega_{Z}^{\bullet}/torsion) \ar[r]^{u} \ar[rd]^{\beta} &  f^{*}(\Omega_{Y}^{\bullet}/torsion) \ar[d]^{v}\\
 (g\circ f)^{*}(\Omega^{\bullet}/torsion) \ar[rr]^{\gamma} & \quad & \Omega_{X}^{\bullet}/torsion } .$$
 Here $\alpha$ is surjective and the image of $\beta$ is the sub-sheaf $f^{**}(g^{**}(\Omega_{Z}))$ by definition. Also the image of $\gamma$ is $(g\circ f)^{**}(\Omega_{Z}^{\bullet})$ by definition. Now the commutativity of the diagram allows to conclude.$\hfill \blacksquare$\\

 \newpage

\section{Definition of the sheaf $\alpha_{X}^{\bullet}$}

The following result is the key of the definition of the sheaf $\alpha_{X}^{\bullet}$ on the reduced  complex space $X$.

\begin{thm}\label{equiv.}
Let $X$ be a reduced complex space and let $S$ be a closed analytic subset with no interior point in $X$ containing the singular set of $X$. Let
 $\alpha$ be a section on $X$ of the sheaf $\omega_{X}^{p}$. The following properties are equivalent for $\alpha$:
\begin{itemize}
\item There exists locally on $X$ a desingularization $\tau : \tilde{X} \to X$ such that $\alpha$ extends to a section on $X$ of  the sub-sheaf $\tau_{*}\tau^{**}(\Omega_{X}^{p})$ of $\omega_{X}^{p}$.\hfill $(A)$
\item There exists locally on $X$ a finite collection $(\rho_{j})_{j \in J}$ of $\mathscr{C}^{\infty}$ functions on $X \setminus S$ which are locally bounded near $S$ and holomorphic $p-$forms $(\omega_{j})_{j \in J}$ on $X$ such that $\alpha = \sum_{j \in J} \ \rho_{j}.\omega_{j}$ as a $(p, 0)$ currents on $X$.\hfill $(B)$
\end{itemize}
\end{thm}

Note that under the second property stated in the theorem, the $(p, 0)-$current on $X$ associated to the form $\sum_{j \in J} \ \rho_{j}.\omega_{j}$ on $X \setminus S$ is defined by
$$ \mathscr{C}_{c}^{\infty}(X)^{n-p, n} \ni \varphi \mapsto \int_{X} \ \varphi\wedge( \sum_{j \in J} \ \rho_{j}.\omega_{j})$$
and this integral is absolutely convergent as the functions $\rho_{j}$ are locally bounded near a point in $S$. It defines a $(p, 0)-$current on $X$ with order $0$. The assumption that $\alpha$ is a section of the sheaf $\omega_{X}^{p}$ implies that this current is $\bar\partial-$closed on $X$.

\parag{Proof} Let us begin by the implication $(A) \Rightarrow (B)$. By definition, a section $\alpha \in \omega_{X}^{p}$ is in the sub-sheaf $\tau_{*}\tau^{**}(\Omega_{X}^{\bullet})$ if,  locally on $\tilde{X}$, it can be written as a linear combination of pull-back of holomorphic forms on $X$ with holomorphic coefficients in $\mathcal{O}_{\tilde{X}}$. Using the properness of the desingularization $\tau$ and a $\mathscr{C}^{\infty}$ partition of the unity on $\tilde{X}$ we obtain the first part of  $(B)$ because $\tau$ induces an isomorphism $\tilde{X} \setminus \tau^{-1}(S) \to X \setminus S$ by hypothesis. The last property in $(B)$, that is to say the fact that the current defined on $X$  by the right hand-side coincides with $\alpha$, is consequence of the fact that both are sections of the sheaf $\omega_{X}^{p}$ and are equal on $X \setminus S$.

 Note that the formula given describes the direct image by $\tau$  of  the $(p, 0)-$current defined by $\alpha$ on $\tilde{X}$.\\
  To prove the implication $(B) \Rightarrow (A)$ consider the pull-back to $\tilde{X} \setminus \tau^{-1}(S)$ of the form $\sum_{j \in J} \ \rho_{j}.\omega_{J}$. We obtain a holomorphic form on $\tilde{X} \setminus \tau^{-1}(S)$  which has locally bounded coefficients along $\tau^{-1}(S)$ when we compute it in a local chart of $\tilde{X}$ near a point of $\tau^{-1}(S)$. So it extends to a holomorphic form  $\tilde{\alpha}$ in $\Omega_{\tilde{X}}^{p}$ and then defines a section $\tau_{*}(\tilde{\alpha})$ of the sheaf $L_{X}^{p}$. But $\tau_{*}(\tilde{\alpha})$ and $\alpha$ coincide on $X \setminus S$, so on $X$ as the sheaf $\omega_{X}^{p}$ has no non zero section supported in $S$.\\
But we have a better estimate. Consider, locally on $\tilde{X}$, a holomorphic function  $g \in \mathcal{O}_{\tilde{X}}$ such that $g.\Omega_{\tilde{X}}^{p} \subset \tau^{**}(\Omega_{X}^{p}) $ (such a non zero $g$ existe as the sheaf $\Omega_{\tilde{X}}^{p}\big/\tau^{**}(\Omega_{X}^{p})$ has its support in $\tau^{-1}(S)$), then the local boundness assumption on the functions $\rho_{j}$ near any point in $S$ gives a local constant $C$ on $X$ such that our section $\tilde{\alpha}$ has coefficients bounded by $C.\vert g\vert$ in a local chart of $\tilde{X}$ on which $\Omega_{\tilde{X}}^{p}$ is a free $\mathcal{O}_{\tilde{X}}-$module.\\
The following lemma allows to conclude that $\tilde{\alpha}$ is a section of the sheaf $\tau^{**}(\Omega_{X}^{p})$ and then $\alpha = \tau_{*}(\tilde{\alpha})$ is a section of  the sheaf $\tau_{*}\tau^{**}(\Omega_{X}^{p})$.$\hfill\blacksquare$

\begin{lemma}\label{cloture}
Let $M$ be a normal complex space and let $\mathcal{F} \subset \mathcal{O}_{M}^{N}$ be a coherent sheaf such that the quotient sheaf $\mathcal{O}_{M}^{N}\big/\mathcal{F} $ has support in a closed analytic subset $S$ without interior point in $M$. Let $\sigma \in \Gamma(M, \mathcal{O}_{M}^{N})$ such that for any local system of generators $g_{1} \dots, g_{k}$ of the annihilator of the sheaf $\mathcal{O}_{M}^{N}/\mathcal{F}$ in a neighbourhood of a compact set $K$ in $M$ there exists a constant $C$ such that we have $\| \sigma \| \leq C.\sup_{i\in [1, k]} \vert g_{i}\vert$ at each point of $K$.
 Then $\sigma$ is in $\Gamma(M, \mathcal{F})$.
\end{lemma}

\parag{Proof} Let $\tau : \tilde{M} \to M$ be a proper modification of $M$ with $\tilde{M}$ normal, such that the sheaf $\tau^{**}(\mathcal{F}) \subset \mathcal{O}_{\tilde{M}}^{N}$ becomes  locally free (of rank $N$). Then we have, locally on $\tilde{M}$ a holomorphic function $\gamma \in \mathcal{O}_{\tilde{M}}$ such that   $\gamma.\mathcal{O}^{N}_{\tilde{M}} = \tau^{**}(\mathcal{F})$. Fix a point $y_{0}$ in $\tilde{M}$ and note $x_{0} := \tau(y_{0})$. Let $g_{1}, \dots, g_{k}$ be holomorphic function in an open neighbourhood of $x_{0}$ which generates the ideal $\mathcal{J}$ in $ \mathcal{O}_{M}$ which annihilates $\mathcal{O}_{M}^{N}\big/\mathcal{F} $. Then, near $y_{0}$ there exists $i \in [1, k]$ and a holomorphic  invertible function $\theta$ near $y_{0}$ such that $\gamma = g_{i}.\theta$ near $y_{0}$. As we have $\| \sigma \| \leq C_{i}.\vert g_{i}\vert $ near $x_{0}$ by hypothesis, there exists some constant $\tilde{C}$ with $\| \tau^{**}(\sigma)\| \leq \tilde{C}.\vert \gamma \vert$ near $y_{0}$. Then the meromorphic section $\tau^{**}(\sigma)/\gamma$ of the sheaf $\mathcal{O}_{M}^{N}$  is holomorphic near $y_{0}$ by normality of $\tilde{M}$ and then $\tau^{**}(\sigma)$ belongs to $\tau^{**}(\mathcal{F})$. But, as $\tau_{*}(\tau^{**})$ induces the identity on the sheaf $\mathcal{O}_{M}$, $\sigma $ belongs to $\Gamma(M, \mathcal{F})$.$\hfill \blacksquare$\\

\begin{defn}\label{inv.}
The graduate sub-sheaf $\alpha^{\bullet}_{X} := \tau_{*}(\tau^{**}(\Omega_{X}^{\bullet}))$ of the sheaf $L_{X}^{\bullet}$ is independent of the choice of the desingularization. So it is naturally defined and $\mathcal{O}_{X}-$coherent thanks to Grauert's direct image theorem.
\end{defn}

\parag{Remark} A direct proof of the fact that the sub-sheaf $\tau_{*}(\tau^{**}(\Omega_{X}^{\bullet}/torsion)$ is independent of the choice of the  desingularization $\tau$ is rather easy. But it does not gives the implication $(B) \Rightarrow (A)$, and property $(B)$ is a simple characterization of the sheaf $\alpha_{X}^{\bullet}$ which does not use any desingularization of $X$.\\
Nevertheless neither $(A)$ nor $(B)$ are easy to compute on examples (see section 6).
\begin{cor}\label{reducible case}
Let $X$ be a pure dimensional reduced complex space and let
  $$X := \cup_{i\in I} X_{i}$$
 be its decomposition in irreducible components. Then the sheaf $\alpha_{X}^{\bullet}$ has a natural injection  in the locally finite  direct sum of the direct images in $X$ of the sheaves $\alpha_{X_{i}}^{\bullet}$ for $i \in I$.
\end{cor}

\parag{Proof} This is an easy consequence of the fact that a section of the sheaf $L_{X}^{\bullet}$ is a section of $\alpha_{X}^{\bullet}$ if and only if it satisfies the condition $(B)$ in the previous theorem, because we have an isomorphism  $L_{X}^{\bullet} \simeq  \oplus_{i \in I} \, ( j_{i})_{*}(L_{X_{i}}^{\bullet})$, where $j_{i} : X_{i}. \to X$ is the inclusion.$\hfill \blacksquare$\\

Note that when $X$ is not irreducible the injective map $\alpha_{X}^{\bullet} \to \oplus_{i \in I} \, (j_{i})_{*}(\alpha_{X_{i}}^{\bullet})$ is not an isomorphism, in general, because the injective map
$$ \Omega_{X}^{\bullet}/torsion  \hookrightarrow \oplus_{i \in I}\, (j_{i})_{*}(\Omega_{X_{i}}^{\bullet}/torsion)$$
is not an isomorphism, in general.\\
Remark that, for each $i \in I$, and any point $x \in X_{i}$, the ``restriction'' map
 $$\alpha_{X,x}^{\bullet} \to \alpha_{X_{i}, x}^{\bullet}$$
is surjective because each restriction map $\Omega_{X, x}^{\bullet}/torsion \to \Omega_{X_{i}, x}^{\bullet}/torsion$ is surjective.\\

\section{Universal pull-back for $\alpha_{X}^{\bullet}$}

The main result of this paragraph if the following theorem.

\begin{thm}\label{pull-back}
For any holomorphic map $f : X \to Y$ between reduced complex spaces, there exists a functorial\footnote{We shall make this precise in the theorem \ref{pull-back bis} below.} graduate $\mathcal{O}_{Y}-$morphism
$$ \hat{f}^{*} : f^{*}\alpha_{Y}^{\bullet} \to \alpha_{X}^{\bullet} $$
which is compatible with the usual pull-back of the sheaf $\Omega^{\bullet}_{Y}/torsion$.\\
For any holomorphic maps $f : X \to Y$ and $g : Y \to Z$ between reduced  complex spaces we have
 \begin{equation*}
 \hat{f}^{*}(\hat{g}^{*}(\alpha)) = \widehat{g\circ f}^{*}(\alpha) \qquad \forall \alpha \in \alpha_{Z}^{\bullet}.\tag{1}
 \end{equation*}
\end{thm}

\bigskip

Let now give a precise formulation of this result. For that purpose let $\mathcal{C}$ be the category of reduced  complex spaces with morphisms all holomorphic maps. We may enrich this category, using the universal pull-back property for the graded sheaf $\Omega_{X}^{\bullet}/torsion$ :\\
Let $\mathcal{C}_{diff}$ be the category whose objects are pairs $(X, \Omega_{X}^{\bullet}/torsion)$ where $X$ is an object in $\mathcal{C}$ and where the morphisms are given by  pairs $(f, f^{*})$ where $f : X \to Y$ is a morphism in $\mathcal{C}$ and $f^{*} : f^{*}(\Omega_{Y}^{\bullet}/torsion) \to \Omega_{X}^{\bullet}/torsion$ is the graded pull-back by $f$ of holomorphic forms modulo torsion (see section 2). Of course the forget-full functor $G_{0}: \mathcal{C}_{diff} \to \mathcal{C}$ obtained by $(X, \Omega_{X}^{\bullet}/torsion) \mapsto X, (f, f^{*}) \mapsto f$ is an equivalence of category.\\

Then the precise  content of the theorem above is the following result.

\begin{thm}\label{pull-back bis}{\bf [Precise formulation]}
There exists a category $\mathcal{C}_{b-diff}$ whose objects are pairs $(X, \alpha_{X}^{\bullet})$ where $X $ is in $\mathcal{C}$ and where the graded coherent sheaf $\alpha_{X}^{\bullet}$ has been defined in section 2 for any object $X$ in $\mathcal{C}$. The morphisms are given by pairs $(f, \hat{f}^{*})$ for each $f : X \to Y$ a morphism in $\mathcal{C}$ where $\hat{f}^{*} : f^{*}(\alpha_{Y}^{\bullet}) \to \alpha_{X}^{\bullet}$ is the graded $\mathcal{O}_{X}-$linear sheaf map defined by $f$. Moreover, the following properties holds:
\begin{enumerate}
\item For each $X \in \mathcal{C}$ we have a graded $\mathcal{O}_{X}-$linear injection
$$ \eta_{X} : \Omega_{X}^{\bullet}/torsion \to \alpha_{X}^{\bullet}.$$
\item For any morphism $f : X \to Y$ in $\mathcal{C}$ we have a commutative diagram of graded $\mathcal{O}_{X}-$linear maps of sheaves
\begin{equation*}
\xymatrix{f^{*}(\Omega_{Y}^{\bullet}/torsion) \ar[r]^{f^{*}} \ar[d]_{f^{*}(\eta_{Y})} & \Omega_{X}^{\bullet}/torsion \ar[d]^{\eta_{X}}\\ f^{*}(\alpha_{Y}^{\bullet}) \ar[r]^{\hat{f}^{*}} & \alpha_{X}^{\bullet}} \tag{2}
\end{equation*}
where $\hat{f}^{*}$ is the graded $\mathcal{O}_{X}-$linear map of coherent sheaves on $X$ associated to the holomorphic map  $f$.
\end{enumerate}
\end{thm}

\bigskip

Of course the interest of this result comes from the fact that the sheaf $\alpha_{X}^{\bullet}$ is, in general, strictly bigger that the sheaf $\Omega_{X}^{\bullet}/torsion$; see section 6. \\
Note also that the sheaf $L_{X}^{\bullet}$ does not have  such a  functorial  pull-back by  holomorphic maps: let $\tau : \tilde{X} \to X$ be a desingularization of $X \in \mathcal{C}$ and let $x \in X$ be a point such that $\tau^{-1}(x)$ has dimension $\geq 1$. Let  $\omega$ be a holomorphic form near $\tau^{-1}(x)$ in $\tilde{X}$ which does not induce a torsion form on an irreducible component $\Gamma$ of $\tau^{-1}(x)$. Then, because the map $\tau_{\vert \Gamma} : \Gamma \to X$ factorizes by the constant map to $\{x\}$ the pull-back of $\omega$ on $\Gamma$ has to be zero. But this map factorizes also by the inclusion of $\Gamma$ in $\tilde{X}$ and $\tau$. As the pull-back by $\tau$ is injective (by definition of $L_{X}^{\bullet}$), this gives a contradiction. Such an example is given in section 6.3.\\

\bigskip

\parag{Preliminaries} Consider the following situation :  let $Z$ be  a connected complex manifold and consider a proper holomorphic map $\pi : Z \to X$ which is surjective on a reduced (irreducible) complex space $X$. Let $q := \dim Z - \dim X$ and let $k$ be the number of connected components of the generic fibre of $\pi$. Assume that  we have a k\"{a}hler form $\omega$ on $Z$. Then, after a suitable normalization of $\omega$,  the smooth $(q, q)-$form $w := \frac{1}{k}.\omega^{\wedge q}$  is $d-$closed and satisfies the condition $\pi_{*}(w) = 1$ as a $d-$closed $(0,0)-$current on $X$. This is consequence of the fact that in the Stein factorization $ \pi_{0} : Z \to Y, \theta : Y \to X$ of $\pi$, the reduced complex space $Y$ is irreducible and this implies that the generic fibres of $\pi_{0}$ are in the same connected component of the space of $q-$cycles in $Z$. So the volume computed by $\omega^{\wedge q}$ of the connected components of the generic fibres of $\pi$ is constant, and we may normalized $\omega$ in order that this volume is equal to $1$. Then the $d-$closed $(0,0)-$current $\pi_{*}(w)$ on $X$ is equal to $1$ on a dense Zariski open set in $X$. This implies our claim.\\

Assume now that the complex manifold $Z$ has finitely many connected components $Z_{1}, \dots, Z_{r}$ such that the restriction of $\pi$ is surjective on each $Z_{j}$ and such that each $Z_{j}$ has a k\"{a}hler form $\omega_{j}$. We can normalize each $\omega_{j}$ in order that the form $w_{j} := \frac{1}{k_{j}}.\omega^{\wedge q_{j}}$ has integral equal to $1/r.k_{j}$ on each connected component of the generic fibres of $\pi_{j}:= \pi_{\vert Z_{j}}$ and  then the smooth form $w := \sum_{j=1}^{r}\  w_{j}$ satisfies again the condition $\pi_{*}(w) = 1$ and for each $j$ any connected component of the generic fibres of $\pi_{j}$ has the same $w-$volume equal to $1/r.k_{j}$.\\

{\bf In this situation we shall say that the smooth for $w$ on $Z$ satisfies the condition $(@)$}.\\

The proof of the theorem \ref{pull-back bis} will use the following proposition.

\begin{prop}\label{main}
Let $X = \cup_{i \in I} \ X_{i}$ be the decomposition of a reduced complex space $X$ as the union of its irreducible components. Let $Z := \cup_{j \in J} \ Z_{j}$ be a disjoint union of connected complex k\"{a}hler  manifolds. Assume that we have a map $\theta : J \to I$ which is surjective and has finite fibres. Let $\pi : Z \to X$ be a proper holomorphic map such that for each $j \in J$ it induces a surjective map 
$$ \pi_{j} : Z_{j} \to X_{\theta(j)} $$
and let $q_{j}:= \dim Z_{j} - \dim X_{\theta(j)}$.  For each $ i \in  J$ let   $w_{j}$ be a smooth $(q_{j}, q_{j})-$form on $Z_{j}$ which is $d-$closed and satisfies the condition $(@)$ relative to the restriction of $\pi$ to $\pi^{-1}(X_{i}) = \cup_{\theta(j) = i} \ Z_{j}$ (see preliminaries above). Note $w := \sum_{j \in J} \ w_{j}$. \\
Let $\beta$ be a section on $Z$ of the sheaf $\pi^{**}(\Omega_{X}^{p})$. Then we have:
\begin{enumerate}
\item The $\bar\partial-$closed $(p, 0)-$current $\pi_{*}(\beta\wedge w)$ on $X$ is independent of the choices of the forms $w_{j}$, assuming that they are $d-$closed and satisfy the condition $(@)$.
\item The section $\pi_{*}(\beta\wedge w)$ on $X$ of the sheaf $\omega_{X}^{p}$ is  a section of the sub-sheaf $\alpha_{X}^{p}$.
\item If there exists a section $\alpha$ of the sheaf $\Omega_{X}^{p}/torsion$ such that $\beta = \pi^{**}(\alpha)$ on $Z$, then $\alpha = \pi_{*}(\beta\wedge w)$ as a section on $X$ of the sheaf $\omega_{X}^{p}$.
\end{enumerate}
\end{prop}

\parag{Remarks} \begin{enumerate}
\item It is enough to prove assertion 1) and 3) of the proposition above for each map $\pi_{j}, j \in J$ because the sheaf $\omega_{X}^{p}$ is the direct sum of the sheaves $\omega_{X_{i}}^{p}, i \in I$ and the restriction of $\beta$ to $Z_{j}$ is a section of the sheaf $\pi_{j}^{**}(\Omega_{X_{\theta(j)}}^{p})$ for each $j \in J$.\\
This is not the case for the assertion 2) of the proposition : the sheaf $\alpha_{X}^{p}$ is a sub-sheaf of the direct sum of the sheaves $\alpha_{X_{i}}^{p}, i \in I$ but, in general, strictly smaller than this direct sum. Note also that the condition on $\beta$ to be a section of the sheaf $\pi^{**}(\Omega_{X}^{p})$ is stronger than the condition on each $\beta_{j}:= \beta_{\vert Z_{j}}, j \in J$ to be a section of the sheaf $\pi_{j}^{**}(\Omega_{X_{\theta(j)}}^{p})$.
\item In general, a section $\beta \in \Gamma(Z, \pi^{**}(\Omega_{X}^{p}))$ is not equal to some  $\pi^{**}(\alpha)$ where $\alpha$ is in $\Gamma(X, \Omega_{X}^{p})$ even in the case where $\pi : Z \to X$ is a desingularization of $X$.$\hfill \square$
\end{enumerate}

\parag{Proof} Thanks to the previous remark, we may assume that $X$ is irreducible to prove assertions 1) and 3) of the proposition.\\
 In the case $q_{j} = 0$ the map $\pi_{j}$ is generically finite and $w_{j}$ is a locally constant function  on $Z$ with a prescribed value on each $Z_{j}$. So there is no choice for $w_{j}$ and the first assertion of the proposition is trivial. As the second assertion is also clear in this case (the sheaf $\omega_{X}^{p}$ has no torsion on $X$ by definition), we shall assume $q_{j} \geq 1$ in the sequel.\\
  The fact that the current $\pi_{*}(\beta\wedge w)$ is $\bar\partial-$closed on $X$ is consequence of the fact that on each $Z_{j}$ the smooth $(p+q_{j}, q_{j})$ form $\beta\wedge w_{j}$ is $\bar\partial-$closed and of the holomorphy of $\pi$. Let $w'$ be a smooth form on $Z$ which is $d-$closed and satisfies the condition $(@)$. We want to show that $\pi_{*}(\beta\wedge (w-w'))$ vanishes as a section of the sheaf $\omega_{X}^{p}$. Let $X'$ be the open and  dense subset of $X$ of smooth points for which the Stein factorization of each $\pi_{j} : Z_{j} \to X$ is a covering of degree $k_{j}$. Remember that, as we assume that $X$ is irreducible here, the set $I $ is reduced to one point and so $J$ is a finite set. On this open set $X'$ it is enough to prove that for each $j \in J$ the current $(\pi_{j})_{*}(\beta\wedge (w_{j}-w'_{j}))$ vanishes. So we can fix $j$ and replace locally $X'$ by one sheet of the corresponding finite covering and make the proof in this case. That is to say that we may assume that $Z$ is  smooth and connected and that $\pi : Z \to X$ has connected fibers on $X$.\\  
  In this case the generic fibres of $\pi$ are irreducible and of dimension  $q$. For any $x \in X'$ there exists an open neighbourhood $V(x)$ of $\pi^{-1}(x)$ which is a deformation retract of $\pi^{-1}(x)$. Then we have an isomorphism $H^{2q}(V(x), \C) \to \C$ which is given by integration on $\pi^{-1}(x)$. But $w$ and $w'$ have the same integral on $\pi^{-1}(x)$ by the property $(@)$. So there exists a $(2q-1)$ smooth form $\theta$ on $V(x)$ such that $d\theta = w - w'$ by de Rham's theorem.\\
   Consider now a small open neighbourhood $U$ of $x$ in $X'$ such that $\pi^{-1}(U) \subset V(x)$. Let $x_{1}, \dots, x_{n}$ be a local coordinate system on $U$. Then the sheaf $\pi^{*}(\Omega_{X}^{p})$ is a free sheaf of $\mathcal{O}_{Z}-$modules on $\pi^{-1}(U)$ with basis $\pi^{*}(dx^{L})$ where $L$ runs in all ordered sub-sets of cardinal $p$ in $[1, n]$. If we write \ $\beta = \sum_{\vert L\vert = p} \ g_{L}.\pi^{*}(dx^{L})$ on $U$ the holomorphic functions $g_{L}$ on $\pi^{-1}(U)$  are constant along the fibres of $\pi$ and so there exists holomorphic functions $f_{L}, \vert L\vert = p$ with $g_{L} = \pi^{*}(f_{L})$  (recall that $U$ is a smooth open set in $X$). This means that there exists a holomorphic $p-$form $\alpha$ on $U$ such that $\beta = \pi^{*}(\alpha)$ on $\pi^{-1}(U)$. \\
   Let $\psi \in \mathscr{C}_{c}^{\infty}(U)^{(n-p, n)}$. By definition of the direct image we have
 $$ \langle \pi_{*} (\beta\wedge d\theta) , \psi \rangle = \int_{\pi^{-1}(U)} \ \beta\wedge d\theta \wedge \pi^{*}(\psi) .$$
 But it follows from the equality $\beta = \pi^{*}(\alpha)$ on $\pi^{-1}(U)$ that the form
  $$\beta\wedge \pi^{*}(\psi) = \pi^{*}(\alpha\wedge \psi)$$
   is $d-$closed as $ \alpha \wedge \psi$ is $d-$closed on $U$ (its degree is $2n$). So by Stokes formula the integral
 $$  \int_{\pi^{-1}(U)} \ \beta\wedge d\theta \wedge \pi^{*}(\psi)  = \pm  \int_{\pi^{-1}(U)} \ d\big(\beta\wedge \theta \wedge \pi^{*}(\psi) \big) $$
 vanishes. This implies that the section $\pi_{*}(\beta\wedge ( w- w'))$ of the sheaf $\omega_{X}^{p}$ vanishes on the open dense subset  $X'$, so everywhere on $X$ as the sheaf $\omega_{X}^{p}$ has no torsion.\\
 The assertion 3) of the proposition is clear, because the equality is obvious at the generic points in $X$.\\
 
 Let now prove the assertion 2). We no longer assume that $I$ has a unique point. \\
  Let $\tau : \tilde{X} \to X$ be a desingularization of $X$, so $\tilde{X}$ is the disjoint union of desingularizations $\tau_{i} : \tilde{X}_{i}  \to X_{i}$ for each $i \in I$, and consider the commutative diagram
$$\xymatrix{ \tilde{X}\times_{X,str} Z  \ar[d]^{\tilde{\pi}} \ar[r]^{\quad \tilde{\tau}} & Z \ar[d]^{\pi} \\ \tilde{X} \ar[r]^{\tau} & X } $$
where $ \tilde{X}\times_{X,str} Z$ is the strict transform, so the union of  irreducible components of $ \tilde{X}\times_{X} Z$ which dominate some  $\tilde{X}_{i}$.\\
 Then the $p-$form $\tilde{\tau}^{**}(\beta)$ gives, for each such component,  a section of the sheaf  $(\tilde{\tau}\circ \pi)^{**}(\Omega_{X}^{p})$ and as the $d-$closed  form $\tilde{\tau}^{*}(w)$ satisfies the condition $(@)$ for the map $\tilde{\pi}$, the $\bar\partial-$closed current $\tilde{\pi}_{*}(\tilde{\tau}^{**}(\beta)\wedge \tilde{\tau}^{*}(w))$ is in fact a $p-$holomorphic form on $\tilde{X}$ thanks to Dolbeault-Grothendieck's lemma. This already proved that $\alpha := \pi_{*}(\beta\wedge w)$ is a section of the sheaf $L_{X}^{p}$, because $\tau^{**}(\pi_{*}(\beta\wedge w)) = \tilde{\pi}_{*}(\tilde{\tau}^{**}(\beta)\wedge\tilde{\tau}^{*}( w)))$ at the generic points of $\tilde{X}$, so everywhere on $\tilde{X}$.\\
 Now the map $\eta : \Omega_{\tilde{X}}^{p} \to \Omega_{\tilde{X}}^{p}$ given by  $\gamma \mapsto \tilde{\pi}_{*}(\tilde{\pi}^{**}(\gamma)\wedge \tilde{\tau}^{*}(w))$ is the identity map, thanks to the assertion 3). So, if $\tilde{\pi}^{**}(\gamma)$ gives a section of the image of  the sub-sheaf  $\tilde{\pi}^{**}\big(\tau^{**}(\Omega_{X}^{p})) $ of the sheaf $\Omega_{\tilde{X}\times_{X,str} Z}^{p}\big/torsion$, $\gamma$  \ will be  a section of  the image of the sub-sheaf  $\tau^{**}(\Omega_{X}^{p})$ because the map $\tilde{\pi}^{*} : \tilde{\pi}^{*}(\Omega_{\tilde{X}}^{p}) \to \Omega_{(\tilde{X}\times_{X,str}Z) }^{p}$ is injective. \\ 
 
 Apply this to $\gamma := \tau^{**}(\alpha) = \tilde{\pi}_{*}\big(\tilde{\tau}^{**}(\beta)\wedge \tilde{\tau}^{*}(w)\big)$ which is a section of $\Omega_{\tilde{X}}^{p}$ as we already proved that $\alpha $ is a section in $L_{X}^{p}$; we obtain that $\tau^{**}(\alpha)$ is a section of the sheaf $\tau^{**}(\Omega_{X}^{p})$ because, as the diagram above commutes,  $\tilde{\tau}^{**}(\beta)$ is a section of the sheaf  \  $\tilde{\tau}^{**}\big(\pi^{**}(\Omega_{X}^{p})\big) = \tilde{\pi}^{**}\big(\tau^{**}(\Omega_{X}^{p}\big)$ thanks to the lemma \ref{pull-back 0}. $\hfill \blacksquare$\\

\parag{Remark}
If $Z$ is not assumed to be smooth in the previous proposition, replacing $Z$ by a  projective desingularization $\sigma : \tilde{Z} \to Z$ (as before, this means that $\tilde{Z}$ is the disjoint union of projective  desingularizations $\sigma_{j} : \tilde{Z}_{j} \to Z_{j}$ for $j \in J$), the proposition applies to the proper map $\sigma\circ \pi$ and to $\tilde{\beta} := \sigma^{*}(\beta)$ which is a section of the sheaf 
$(\sigma\circ\pi)^{**}(\Omega_{X}^{p})$. Then the result is still true.$\hfill \square$

 \parag{Proof of the theorem \ref{pull-back}}The first step in proving the theorem will be the construction of $\hat{f}^{*}(\alpha) \in \alpha_{X}^{\bullet} $ when $\alpha$ is a section of the sheaf $\alpha_{Y}^{\bullet}$. So let $\alpha$ be a section on $Y$ of the sheaf $\alpha_{Y}^{p}$. Let $\tau : \tilde{Y} \to Y$ be a locally projective desingularization of $Y$. Consider the following  commutative diagram
$$ \xymatrix{\tilde{X}\ar[d]^{\pi} \ar[r]^{\tilde{f}} & \tilde{Y} \ar[d]^{\tau} \\ X \ar[r]^{f} & Y} $$
where $\tilde{X} \subset X\times_{Y}\tilde{Y} $ is the strict transform of $X$, that is to say the union of irreducible components of $X\times_{Y}\tilde{Y} $  which dominate an irreducible component of $X$,  and where $\pi$ and $\tilde{f}$ are induced by the natural projections of $X\times_{Y}\tilde{Y}$.\\
Now the problem is local on $X$ and $Y$ and we may assume that $X, Y$ and $\tilde{Y}$ are k\"{a}hler. So we may assume that we have on $\tilde{X}$ a smooth $d-$closed form $w$ which satisfies the condition $(@)$ for the proper map $\pi$ (we may replace $\tilde{X}$ by a  projective desingularization to reach the precise situation of the proposition \ref{main}; see the remark above).\\
Then, as $\beta $ is a section of $\tau^{**}(\Omega_{Y}^{p})$, the form $\tilde{f}^{**}(\beta)$ is a section of $\pi^{**}(\Omega_{X}^{p})$ because if we write locally on $\tilde{Y}$ 
$$ \beta := \sum_{l} \ g_{l}.\tau^{**}(\omega_{l}) $$
where $\omega_{l}$ are local sections of $\Omega_{Y}^{p}$ and $g_{l}$ are holomorphic functions on $\tilde{Y}$, we obtain
$$ \tilde{f}^{**}(\beta) = \sum_{l} \ \tilde{f}^{*}(g_{l}).\tilde{f}^{**}(\tau^{**}(\omega_{l})) $$
and the equality $\tilde{f}^{**}(\tau^{**}(\omega_{l})) = \pi^{**}(f^{**}(\omega_{l}))$ due to the commutativity of the diagram and the lemma \ref{pull-back 0} shows that $\tilde{f}^{**}(\beta)$ is a section of the sheaf $\pi^{**}(\Omega_{X}^{p})$. So we can apply the proposition \ref{main} and obtain that $\pi_{*}(\tilde{f}^{**}(\beta)\wedge w)$ is a section of the sheaf $\alpha_{X}^{p}$. This will give the definition of $\hat{f}^{*}(\alpha)$ when we will have proved that it is independent of the choice of the desingularization $\tau : \tilde{Y} \to Y$.\\
 Note that the proposition \ref{main} already gives the independence of the choice of $w$ (assumed $d-$closed and satisfying $(@)$) in this construction.\\
 Note also that the proposition \ref{main} gives also that for $\alpha$ a section of $\Omega_{Y}^p/torsion$ $\hat{f}^{*}(\alpha)$ is  a section of $\Omega_{X}^{p}/torsion$ and coincides with the usual pull-back $f^{*}(\alpha)$ (see section 2).\\
 
 Remark now that, as the sheaf $\alpha_{X}^{p}$ has no torsion on $X$, to prove the independence of $\hat{f}^{*}(\alpha)$ on the choice of the desingularization $\tau$, it is enough to prove it at the generic points of $X$. Moreover, this problem is local on $X$ and so we may assume  that $X$ is smooth and connected.\\
 In our construction, we sum the various direct images $(\pi_{j})_{*}(\tilde{f}^{*}(\beta)\wedge w_{j})$ when $j$ describes the various connected components of the desingularization of $\tilde{X}$. Each such component is send by $\tilde{f}$ in a connected component of $\tilde{Y}$ and then it is enough to show the invariance of the current $(\pi_{j})_{*}(\tilde{f}^{*}(\beta)\wedge w_{j})$ if we change only one connected component of $\tilde{Y}$ in the given desingularization, and also if we consider only the corresponding  connected components of the desingularization of $\tilde{X}$. So, in fact, it is enough to prove the following special case of our problem:\\
 Assume that $X$ is smooth and connected and that $Y$ is irreducible. Let $\tau : \tilde{Y} \to Y$ is a projective desingularization of $Y$ and that $\theta : \tilde{\tilde{Y}} \to \tilde{Y}$ is a proper smooth modification of $\tilde{Y}$. So our new desingularization of $Y$ is $\theta\circ \tau :  \tilde{\tilde{Y}} \to Y$. Note that we can make this assumption without any lost of generality because two desingularizations can always be dominated by a third one (thanks to Hironaka).\\
 Now we shall consider the following diagram, where $\tilde{X}$ is a  desingularization of an irreducible component of the strict transform $X\times_{Y}\tilde{Y}$ and $\tilde{\tilde{X}}$ is a desingularization of the strict transform of $\tilde{X}\times_{\tilde{Y}}\tilde{\tilde{Y}}$:
 
 $$ \xymatrix{\tilde{\tilde{X}} \ar[d]^{\tilde{\theta}} \ar[r]^{\tilde{\tilde{f}}}& \tilde{\tilde{Y}} \ar[d]^{\theta}\\ \tilde{X}\ar[d]^{\tilde{\tau}} \ar[r]^{\tilde{f}} & \tilde{Y} \ar[d]^{\tau} \\ X \ar[r]^{f} & Y} $$
 
 Let $q$ the dimension of the generic fibres of $\tilde{\tau}$ and $k$ the number of connected components of its generic fibres. Let $\omega$ be a k\"{a}hler form of $\tilde{Y}$ normalized in order that the form $\tilde{f}^{*}(\omega^{\wedge q})$ satisfies the condition $(@)$ for the map $\tilde{\tau}$. Let $\tilde{q}$ be the dimension of the generic fibre of $\tilde{\theta}$ and let $\tilde{\omega}$ a k\"{a}hler form on $\tilde{\tilde{Y}}$ normalized in order that the form $\tilde{\tilde{f}}^{*}(\tilde{\omega}^{\wedge \tilde{q}})$ satisfies the condition $(@)$ for the map $\tilde{\theta}$. Now consider the $(q+\tilde{q}, q+\tilde{q})-$smooth form $w := \tilde{\tilde{f}}^{*}(\theta^{*}(\omega^{\wedge q})\wedge \tilde{\omega}^{\tilde{q}})$  on $\tilde{\tilde{X}}$ which is $d-$closed. It satisfies the condition $(@)$ for the map $\tilde{\theta}\circ\tilde{\tau}$.\\
 So the definition of $\hat{f}^{*}(\alpha)$ using the desingularization $\theta\circ\tau$ is given by 
 $$(\tilde{\theta}\circ\tilde{\tau})_{*}\big((\tilde{\tilde{f}}\circ \theta)^{*}(\beta)\wedge w\big),$$
 But, as $\tilde{f}^{**}(\beta)$ is a section of the sheaf $\Omega^{p}_{\tilde{X}}/torsion$, we have the equality  
$$\tilde{\theta}_{*}(\tilde{\theta}^{**}(\tilde{f}^{**}(\beta))\wedge\tilde{\tilde{f}}^{*}(\tilde{\omega}^{\tilde{q}})) = \tilde{f}^{**}(\beta)$$ and the conclusion follows from the fact that 
$$( \tilde{\theta}\circ\tilde{\tau})_{*}\big((\tilde{\tilde{f}}\circ\theta)^{**}(\beta)\wedge w\big) = \tilde{\tau}_{*}[\tilde{\theta}_{*}\big(\tilde{\theta}^{**}(\tilde{f}^{**}(\beta))\wedge\tilde{\tilde{f}}^{*}(\tilde{\omega}^{\tilde{q}})\big)\wedge \tilde{f}^{*}(\omega^{q}))].$$

 The compatibility of this construction with the pull-back of holomorphic forms modulo torsion which is given by the last assertion of the proposition \ref{main} obviously gives that the injective $\mathcal{O}_{X}-$linear morphism 
 $$\eta_{X} : \Omega_{X}^{\bullet}/torsion \to \alpha_{X}^{\bullet} $$
  for each $X \in \mathcal{C}$  gives  the commutative diagram $(2)$ of the precise formulation \ref{pull-back bis} of the theorem for each morphism $f : X \to Y$ in $\mathcal{C}$.    \\
  
 Now we have to prove the functoriality of $\hat{f}^{*}$. Then consider a holomorphic maps $f: X \to Y$ and  $g : Y \to Z$. We want to prove the formula $(1)$ of the theorem.\\
 Consider the commutative diagram
 
 $$\xymatrix{\tilde{\tilde{X}} \ar[r]^{\tilde{\tilde{f}}} \ar[d]^{\tilde{\theta}} & \tilde{\tilde{Y}} \ar[d]^{\theta} &\quad\\
 \tilde{X} \ar[r]^{\tilde{f}} \ar[d]^{\tau_{2}}& \tilde{Y} \ar[r]^{\tilde{g}} \ar[d]^{\tau_{1}} & \tilde{Z} \ar[d]^{\tau} \\
 X \ar[r]^{f} & Y \ar[r]^{g} & Z } $$
 where $\tau : \tilde{Z} \to Z$ is a desingularization, where $\tilde{g} : \tilde{Y} \to \tilde{Z}$ is the strict transform of $g$ by $\tau$, where $\tilde{f} : \tilde{X} \to \tilde{Y}$ is the strict transform of $f$ by $\tau_{1}$, where $\tilde{\tilde{f}} : \tilde{\tilde{X}} \to \tilde{\tilde{Y}}$ is the strict transform of $\tilde{f}$ by $\theta : \tilde{\tilde{Y}} \to \tilde{Y}$ which is a desingularization of $\tilde{Y}$.\\
 Let $\alpha$ be a section of $\alpha_{Z}^{p}$, note $\beta := \tau^{**}(\alpha) \in \tau^{**}(\Omega_{Z}^{p})$\footnote{See the simple lemma \ref{satur} below.} and let $w_{1}$ and $w_{2}$ be smooth $d-$closed forms satisfying the condition $(@)$ of the proposition \ref{main} for the maps $\tau_{1}$ and $\tilde{\theta}$ respectively. We have 
 $$ \hat{g}^{*}(\beta) = (\tau_{1})_{*}(\tilde{g}^{**}(\beta) \wedge w_{1})$$
 but we have also, because $\tilde{g}^{**}(\beta) $ is a section of $\tau_{1}^{**}(\Omega_{Y}^{p})$
 $$ \hat{g}^{*}(\beta) = (\theta\circ\tau_{1})_{*}(\theta^{**}(\tilde{g}^{**}(\beta) \wedge \theta^{*}(w_{1}))).$$
 Then we obtain
 $$ \hat{f}^{*}(\hat{g}^{*}(\alpha)) = (\tilde{\theta}\circ \tau_{2})_{*}\big(\tilde{\tilde{f}}^{**}(\theta^{**}(\tilde{g}^{**}(\beta)))\wedge \tilde{\tilde{f}}^{*}(w_{1})\wedge w_{2}\big) .$$
 As the square
 $$ \xymatrix{\tilde{X} \ar[r]^{\tilde{f}\circ\tilde{g}} \ar[d]_{\tau_{2}} & \tilde{Z} \ar[d]^{\tau}\\
 X \ar[r]^{f\circ g} & Z} $$
 is also the strict transform of $f\circ g$ by $\tau$ we have
 $$ \widehat{f\circ g}^{*}(\alpha) = (\tau_{2})_{*}\big((\tilde{f}\circ\tilde{g})^{**}(\beta)\wedge \tilde{f}^{*}(w_{1})\big) .$$
 Then the conclusion follows from the equality
 $$ \tilde{\theta}_{*}\big(\tilde{\tilde{f}}^{**}(\theta^{**}(\tilde{g}^{**}(\beta)))\wedge \tilde{\tilde{f}}^{**}(\theta^{*}(w_{1}))\wedge w_{2}\big) = \tilde{f}^{**}(\tilde{g}^{**}(\beta))\wedge \tilde{f}^{*}(w_{1})$$
 obtained by the comparaison of both hand-sides at the generic points of $\tilde{X}$.$\hfill \blacksquare$\\

Our next result shows that the sheaf $\alpha_{X}^{\bullet}$ is ``maximal'' in order to construct the pull-back via the method of the proposition \ref{main}.

\begin{thm}\label{pull-back 2}
Let $\pi : Z \to X$ be a proper surjective holomorphic map between irreducible complex spaces. Put $q := \dim Z - \dim X$. Let $\beta \in \alpha_{Z}^{p}$ be  equal to $\hat{\pi}^{*}(\alpha)$ for a  section $\alpha$  of the sheaf $\alpha_{X}^{p}$. Let also $w$ be a smooth $(q, q)-$form on $Z$ which is $d-$closed and satisfies the condition $(@)$ of the proposition \ref{main}. Then the $(p, 0)-$current $\pi_{*}(\beta\wedge w)$ on $X$ (which is $\bar\partial-$closed  and  independent of the choice of $w$ satisfying  $dw = 0$ and $(@)$; see proposition \ref{main})  is equal to the image in $\omega_{X}^{p}$  of the section $\alpha$ of the sheaf $\alpha_{X}^{p}$.
\end{thm}

Using the ``pull-back'' theorem \ref{pull-back} the theorem above will follow from the following simple lemma.

\begin{lemma}\label{satur}
Let $X$ be a complex space and $\tau : \tilde{X} \to X$ a desingularization of $X$. Then the image of the pull-back $\hat{\tau}^{*} : \tau^{*}(\alpha_{X}^{\bullet}) \to \alpha_{\tilde{X}}^{\bullet} = \Omega_{\tilde{X}}^{\bullet}$ is the subsheaf $\tau^{**}(\Omega_{X}^{\bullet})$ of $ \Omega_{\tilde{X}}^{\bullet}$.
\end{lemma}

\parag{Proof} By definition, a section of  this image is locally on $\tilde{X}$ a $\mathcal{O}_{\tilde{X}}-$linear combination of holomorphic forms on $\tilde{X}$ which are locally $\mathcal{O}_{\tilde{X}}-$linear combinations of pull-back by $\tau$ of holomorphic forms on $X$. So the conclusion is clear.$\hfill \blacksquare$

\parag{Proof of the theorem \ref{pull-back 2}} Let $\tilde{\pi} : \tilde{Z} \to \tilde{X}$ be the strict transform of $\pi$ by $\tau$, and denote by $\tilde{\tau} : \tilde{Z} \to Z$ the corresponding projection on $Z$ which is a modification. So we have the following commutative diagram
$$\xymatrix{ \tilde{Z} \ar[r]^{\tilde{\pi}} \ar[d]^{\tilde{\tau}} & \tilde{X} \ar[d]^{\tau} \\ Z \ar[r]^{\pi} & X} $$

 The $(q, q)-$form $\tilde{\tau}^{*}(w)$ is smooth and $d-$closed in $\tilde{Z}$ and satisfies the condition $(@)$ of the proposition \ref{main} for the proper surjective holomorphic map $\tau\circ \tilde{\pi}$.
As we can write  $\beta = \hat{\pi}^{*}(\alpha)$  where  $\alpha$ is a section of $\alpha_{X}^{p}$, we have, by functoriality of the pull-back for the sheaf $\alpha_{Y}^{\bullet}$ and the equality $\tau\circ \tilde{\pi} = \pi\circ\tilde{\tau}$
$$ \hat{\tau}^{*}(\beta) = \widehat{\tilde{\pi}}^{*}(\hat{\tau}^{*}(\alpha)) .$$
But, thanks to the previous lemma, we have $\hat{\tau}^{*}(\alpha) $ which is a section of $\tau^{**}(\Omega_{X}^{p})$ and using the smoothness of $\tilde{X}$ we have $ \widehat{\tilde{\pi}}^{*} = \tilde{\pi}^{**}$. Then we obtain, using the lemma \ref{pull-back 0}, the fact that $\hat{\tau}^{*}(\beta)$ is a section of the sheaf $(\tau\circ\tilde{\pi})^{**}(\Omega_{X}^{p})$. Then the proposition \ref{main}  applies to the map $\tau\circ\tilde{\pi} : \tilde{Z} \to X$ with the form $\tilde{\tau}^{*}(w)$ and the section $\hat{\tau}^{*}(\beta)$ of the sheaf $(\tau\circ\tilde{\pi})^{**}(\Omega_{X}^{p})$ and gives that the $(p, 0)-$current on $X$ given by $\sigma := (\tau\circ\tilde{\pi})_{*}(\hat{\tau}^{*}(\beta)\wedge \tilde{\tau}^{*}(w))$ is $\bar\partial-$closed on $X$ and is a section of the sheaf $\alpha_{X}^{p}$.\\
But the $(p+q, q)-$current $\tilde{\tau}_{*}(\hat{\tau}^{*}(\beta)\wedge \tilde{\tau}^{*}(w))$ is equal to $\beta\wedge w$ at least over the generic points in $X$, the $(p, 0)-$current  $\bar\partial-$closed in $X$  $\pi_{*}(\beta\wedge w)$ is  generically equal to $\sigma$ and $\alpha$. So $\alpha $ and $\sigma$ are equal as sections of the sheaf $\alpha_{X}^{p}$.$\hfill \blacksquare$\\

\section{Integration on cycles}

 \parag{Notations} Let $V$ be a complex manifold and $h$ be a continuous  hermitian form on $V$. So $h$ is a real continuous definite positive $(1,1)-$ differential form on $V$. If $\omega$ is a continuous  $(p,p)-$form on $V$, we shall consider $\omega$ as a continuous sesqui-linear form on $\Lambda^{p}(T_{V})$ and we shall write
 $$ \| \omega \|_{K} \leq C.h^{\wedge p}  $$
 where $K$ is a subset in $V$ and $C > 0$ a constant, if for any point $x \in K$ and any $v_{1}, \dots, v_{p} \in T_{V,x}$ the inequality
 $$ \vert \omega(x)[v_{1}\wedge \dots \wedge v_{p}] \vert \leq C.h^{\wedge p}(x)[v_{1}\wedge \dots \wedge v_{p}]  $$
 holds. For instance, if $\alpha, \beta  \in \Omega_{V}^{p}$ we shall write $\|\alpha \wedge \bar\beta \|_{K} \leq C_{K}.h^{\wedge p}$ when for any $x \in K$ and  any $v_{1}, \dots, v_{p} \in T_{V,x}$ we have
 \begin{equation*}
  \vert \alpha(x)[v_{1}\wedge \dots \wedge v_{p}]\vert.\vert \beta(x)[v_{1}\wedge \dots \wedge v_{p}]\vert \leq C_{K}.h^{\wedge p}(x)[v_{1}\wedge \dots \wedge v_{p}]. \tag{1}
  \end{equation*}
  
  \parag{Remark}
 If $f : W \to V$ is a  holomorphic map and if  $(1)$ holds  then we shall have 
 \begin{equation*}
  \| f^{*}(\alpha)\wedge \overline{f^{*}(\beta)}\|_{f^{-1}(K)} \leq C_{K}.f^{*}(h)^{\wedge p} \tag{2}
  \end{equation*}
 but, in general, $f^{*}(h)$ is no longer definite positive on $W$. \\
 Conversely if $(2)$ holds on a set  $L$ in $W$ then $(1)$ is satisfied on $f(L)$.
 
 \begin{prop}\label{bound}
 Let $X$ be a  reduced complex space, let $S$ be the singular set in $X$  and let $h$ be a continuous hermitian metric on $X$. Let $U$ be a relatively compact open set in $X$. For all $\alpha, \beta \in \alpha_{X}^{p}$ there exists a constant $C_{U} > 0$ such that the following inequality holds at each point in $\bar U \setminus S$
 $$ \| \alpha\wedge \bar\beta \|_{\bar U \setminus S} \leq C_{U}.h^{\wedge p}_{\bar U\setminus S}.$$
 \end{prop}
 
 \parag{Proof}  Remark that the problem is local on the compact set $\bar U \cap S$ because near smooth points in $X$ the assertion obviously holds. Let $\tau : \tilde{X} \to X$ be a desingularization of $X$. Then we shall show that for each point $y \in \tau^{-1}(\bar U \cap S)$ there exists an open neighbourghood $W$ of $y$ in $\tilde{X}$ and a positive constant $C_{W}$ such that the inequality 
 $$\| \tau^{**}(\alpha)\wedge \overline{ \tau^{**}(\beta)} \|_{W} \leq C_{W}.\tau^{**}(h)^{\wedge p}$$
 holds : if $y$ is a point in $\tilde{X}$ we can write in an open neighbourghood $W$ of $y$ 
 $$\alpha = \sum_{\vert I\vert = p} \ g_{I}.\tau^{**}(dx^{I}) \quad {\rm and} \quad  \beta = = \sum_{\vert I\vert = p} \ h_{I}.\tau^{**}(dx^{I}) $$
 where $x_{1}, \dots, x_{N}$ are local coordinates in a closed embedding of an open set $U \subset\subset X$ in $\C^{N}$ near $\tau(y)$. Our estimates is consequence of the facts that the holomorphic functions $g_{I}$ and $h_{I}$ are locally bounded and that for any $(I, J)$ there is a constant $c_{U}^{I, J} > 0$ with
  $$\| dx^{I}\wedge \overline{dx^{J}})\|_{U}  \leq   c_{U}^{I, J}.h^{\wedge p}$$
  because we can assume that $h$ in induced by a continuous hermitain form on $\C^{N}$.\\
  Now the properness of $\tau$ allows to find a  a constant $C_{U}$ such that the inequality 
  $$\| \tau^{**}(\alpha)\wedge \overline{ \tau^{**}(\beta)} \|_{K} \leq C_{U}.\tau^{**}(h)^{\wedge p} $$
  holds on the compact set $K := \tau^{-1}(\bar U)$. This allows to conclude thanks to the remark above.$\hfill \blacksquare$\\
  
  \begin{cor}\label{int 1}
  Let $X$ be a complex space of pure dimension $n$,  and let $\alpha, \beta$ be sections on $X$ of the sheaf $L_{X}^{n}$. Then, if $\rho$ is a continuous compactly supported function on $X$ the integral
  $$ \int_{X\setminus S} \ \rho.\alpha\wedge \bar\beta $$
  is absolutely convergent for any closed analytic subset $S$ containing the singular set in $X$ and its value does not depends on the choice of $S$.\\
  If,  moreover, $\alpha$ and $ \beta$ are sections of the sheaf $\alpha_{X}^{n}$,  for any continuous hermitian metric $h$ on $X$ there is constant $C > 0$ depending on $\alpha, \beta$ and on the support  $K$ of $\rho$ such that 
  \begin{equation*}
   \vert  \int_{X\setminus S } \ \rho.\alpha\wedge \bar\beta \vert \ \leq \  C.\int_{K} \ \vert \rho \vert.h^{\wedge n} \leq C.\vert\vert \rho \vert\vert. \int_{K} \ h^{\wedge n} . \tag{3}
   \end{equation*}
  \end{cor}
  
  \parag{Proof} The first part is consequence of the fact that $\tau^{**}(\alpha)$ and $\tau^{**}(\beta)$ are holomorphic $n-$forms on $\tilde{X}$. The estimates when $\alpha, \beta$ are sections of $\alpha_{X}^{n}$  is a direct consequence of the previous proposition.$\hfill \blacksquare$\\
  
  \begin{defn}\label{int 2}
  For $\alpha, \beta$ sections of the sheaf $L_{X}^{n}$ the common values of the absolutely convergent integrals  $ \int_{X \setminus S} \ \rho.\alpha\wedge \bar\beta$ will be denoted simply  by $ \int_{X} \ \rho.\alpha\wedge \bar\beta$.
  \end{defn}
  
  \begin{lemma}\label{change var.}
  Let $f : Y \to X$ a proper generically finite and surjective holomorphic map between two complex spaces of pure dimension $n$; let $k$ be the generic degree of $\pi$. Let $\alpha, \beta$ be sections on $X$ of the sheaf $L_{X}^{n}$ and $\rho \in \mathscr{C}_{c}^{0}(X)$. Then the holomorphic $n-$forms $f^{**}(\alpha)$ and $f^{**}(\beta)$ are well defined on a dense Zariski open set in $Y$ and extend as sections on $Y$ of the sheaf $L_{Y}^{n}$. We have the equality
  $$ \int_{X} \ \rho.\alpha\wedge \bar\beta = k.\int_{Y} \ f^{*}(\rho). f^{**}(\alpha)\wedge \overline{f^{**}(\beta)} .$$
  \end{lemma}
  
  \parag{Proof} Remark that it is enough to prove the lemma for $\alpha = \beta$. As $\tau^{**}(\alpha)$ is an holomorphic $n-$form on $\tilde{X}$ implies that $\alpha$ is locally $L^{2}$ on $X$. Let $S_{\varepsilon}$ be an  open $\varepsilon-$neighbourhood of $S$ a closed analytic subset in $X$ such that the map $f : Y \setminus f^{-1}(S) \to X\setminus S$ is a finite covering between two complex manifolds. Then the usual change of variable gives, if $\rho$ is in $ \mathscr{C}_{c}^{0}(X)$
  $$ \int_{X \setminus S_{\varepsilon}} \ \rho.\alpha\wedge \overline{\alpha} =  k.\int_{Y \setminus f^{-1}(S_{\varepsilon})} \ f^{*}(\rho).f^{**}(\alpha) \wedge \overline{f^{**}(\alpha)}.$$
  Letting $\varepsilon $ goes to $0$  shows that $f^{**}(\alpha)$ is locally $L^{2}$ on any desingularization of $Y$ and so $f^{**}(\alpha)$  is a section of the sheaf $L_{Y}^{n}$. The conclusion follows easily.$\hfill \blacksquare$\\
  
  \begin{defn}\label{int. cycle}
  Let $X$ be a complex space and let $Y \subset X$ be an irreducible $p-$dimensional analytic subset in $X$. We shall denote $j : Y \to X$ the  the inclusion. Let $\alpha, \beta$ be sections of the sheaf $\alpha_{X}^{p}$ on $X$ and $\rho$ be a continuous function with  compact support in $X$. We shall define the number $\int_{Y} \ \rho.\alpha\wedge \overline{\beta} $ as the integral 
  $$\int_{Y} \ j^{*}(\rho).\hat{j}^{*}(\alpha)\wedge \overline{\hat{j}^{*}(\beta)} .$$
  \end{defn}
  
  Note that this definition makes sense because the pull-back $\hat{j}^{*} : j^{*}(\alpha_{X}^{p}) \to \alpha_{Y}^{p}$ is well defined and because the inclusion $\alpha_{Y}^{p} \subset L_{Y}^{p}$ allows to use the definition \ref{int 2}. Remark that this definition only depends on the irreducible analytic subset $Y$ of $X$. So we may extend by additivity the definition of the integral
  $$ \int_{Y} \ \rho.\alpha\wedge \overline{\beta} $$
   to any $p-$dimensional cycle $Y$ in $X$.

   \parag{Remark} It is not clear that the definition \ref{int. cycle} may be extended to sections of the sheaf $L_{X}^{p}$ for $0 < p < n$, for a $p-$cycle $Y$ contained in the singular locus of $X$.\\ 

 The next lemma shows that the change of variable holds for such a integral.
  
  \begin{lemma}\label{chgt var.}
  Let $f : X \to Y$ be a holomorphic map and let $\alpha, \beta$ be sections on $Y$ of the sheaf $\alpha_{Y}^{p}$. Let $\rho$ be a continuous compactly supported function on $Y$. Let $Z$ be a $p-$cycle in $X$ and assume that the cycle $f_{*}(Z)$ is defined in $Y$\footnote{This means that the restriction of $f$ to $\vert Z\vert$ is proper; see [B-M 1] chapter IV.}. Then the restriction to $\vert Z\vert$ of the continuous function $f^{*}(\rho)$ has compact support and the integral $\int_{Z} \ f^{*}(\rho).\hat{f}^{*}(\alpha)\wedge \overline{\hat{f}^{*}(\beta)} $ is well defined and we have
  $$ \int_{Z} \ f^{*}(\rho).\hat{f}^{*}(\alpha)\wedge \overline{\hat{f}^{*}(\beta)} = \int_{f_{*}(Z)} \ \rho.\alpha\wedge \bar\beta .$$
  \end{lemma}
  
  \parag{proof} First remark that any irreducible component $\Gamma$ of $Z$ which has an image of dimension at most equal to $p-1$ does not contribute to the right hand-side and also to the left hand-side because the forms $\hat{f}^{*}(\alpha)$ and $\hat{f}^{*}(\beta)$ vanish on such a irreducible component :  \\
  Let $g : \Gamma \to f(\Gamma)$ be the map induced by $f$; by functoriality of the pull-back $\hat{g}^{*}$ factorizes through $\alpha_{f(\Gamma)}^{p}$ which is zero.\\
  Then the result is in fact a local statement near each point of the set $\vert f_{*}(Z)\vert$. And because of our previous remark and the fact that closed analytic subsets with no interior point can be neglected in the integrals, it is enough to prove the result when $Z$ is smooth and when $f$ induces an isomorphism of $Z$ on $f(Z)$. In this case, which is not trivial because $Z$ and $f(Z)$ can be contained in the singular sets of $X$ and $Y$, the functorial property of the pull-back and the fact that for a complex manifold $V$ we have $\alpha_{V}^{p} = \Omega_{V}^{p}$ allow to conclude.$\hfill \blacksquare$\\

  \begin{thm}\label{int. param.}
  Let $X$ be a  complex space and let  $(Y_{t})_{t \in T}$ be an analytic family of $p-$cycles in $X$ parametrized by a reduced complex space $T$. Let $\alpha, \beta$ be sections of the sheaf $\alpha_{X}^{p}$ on $X$ and $\rho$ be a continuous function with  compact support in $X$. Then the function $\varphi : T \to \C$ defined by
  $$ \varphi(t) := \int_{Y_{t}} \ \rho.\alpha\wedge \overline{\beta} $$
  is locally bounded and for any given hermitian metric $h$ on $X$  there exists a constant $C$ such the following estimate holds for each $t \in T$:
  \begin{equation*}
  \vert \varphi(t) \vert \leq C.\int_{Y_{t}} \vert \rho \vert.h^{\wedge p} \tag{E}.
  \end{equation*}
  Moreover for each point $t_{0} \in T$ there exists an open neighbourhood $T_{0}$ of $t_{0}$ in $T$ and a closed analytic subset $\Theta_{0}$ with no interior point in $T_{0}$ such that $\varphi$ is continuous on $T_{0}\setminus \Theta_{0}$.
   \end{thm}
  
  \parag{Proof} We shall cut this proof in several steps.
  \parag{Step 1} Let $\nu : \tilde{T}  \to T$ the normalization of $T$. The family $(Y_{\nu(\tilde{t})})_{\tilde{t} \in \tilde{T}}$ is an analytic family of $p-$cycles in $X$ parametrized by $\tilde{T}$, and if the theorem is proved for this family it implies the result for the initial family, because the function is constant on the fibres of the normalization map.\\
  So we shall assume that $T$ is normal in the sequel.
  \parag{Step 2} If the generic cycle $Y_{t}$ is not reduced and irreducible, the normality of $T$ allows to write the family $(Y_{t})_{t\in T}$ as a finite sum of analytic families of $p-$cycles in $X$ parametrized by $T$ such the sum of these families is our initial family (see [B-M 1] ch. IV theorem 3.4.1). So it is enough to prove the theorem for such a family.\\
  We shall assume that for $t$ generic in $T$ the cycle $Y_{t}$ is reduced and irreducible.
  \parag{Step 3} Let $G \subset T \times X$ the cycle-graph of our analytic family. It is a reduced cycle and the projection $\pi : G \to T$ is (by definition) a geometrically flat map, that is to say that there exists an analytic family of cycles $(Z_{t})_{t \in T}$  in $G$ such that for each $t \in T$ we have $\vert Z_{t}\vert = \pi^{-1}(t)$ and such that the generic cycle $Z_{t}$ is reduced and irreducible. Of course, here we have $Z_{t}:= \{t\}\times Y_{t}$ for each $t \in T$.\\
   Note $pr: G \to X$ the projection and define on $G$ the sections of the sheaf $\alpha_{G}^{p}$ by letting $\alpha_{1}:= \hat{pr}^{*}(\alpha)$ and $\beta_{1} := \hat{pr}^{*}(\beta)$. Then, it is enough to prove the theorem for the function $t \mapsto \int_{Z_{t}}\ \tilde{\rho}.\alpha_{1}\wedge \bar{\beta}_{1}$ where $\tilde{\rho} := pr^{*}(\rho)$ thanks to the change variable theorem proved in lemma \ref{chgt var.}. Remark that $pr$ induces an isomorphism of $\vert Z_{t}\vert$ onto $\vert Y_{t}\vert$ for each $t \in T$. Note also that the continuous function $\tilde{\rho}$ on $G$ has a  $\pi-$proper support.
   
\parag{Step 4} Let $\tau : \tilde{G} \to G$ be a desingularization of $G$. Define the subset 
$$\Theta := \{ t \in T\ / \ \exists y \in K \  \dim_{y}\, (\tau\circ pr)^{-1}(t) \geq p+1 \}.$$
This is a locally closed  analytic subset\footnote{See, for instance, the lemma 2.1.8 in [B.15].} in $T$ with no interior point. For $t_{0} \in T$  fix an open neighbourhood $T_{0}$ of $t_{0}$  small enough in order that $\Theta_{0}:= \Theta \cap T_{0}$ is a closed analytic subset. The map 
$$q : \tilde{G}\cap (\tau\circ pr)^{-1}(T_{0})  \setminus (\tau\circ pr)^{-1}(\Theta) \to T_{0} \setminus \Theta_{0}$$
 is $p-$equidimensional on a normal basis, so it is geometrically flat and we have an analytic family $(\tilde{Z}_{t})_{t \in T_{0} \setminus \Theta_{0}}$ of fibres of $q$ which are $p-$cycles in $\tilde{G}$, and for $t$ generic in $T_{0}\setminus \Theta_{0}$ the cycle $\tilde{Z}_{t}$ is irreducible. \\
 Note that the pull-back of $\alpha_{1}$ and $\beta_{1}$ on $\tilde{G} \cap (\tau\circ pr)^{-1}(T_{0})$  are holomorphic $p-$forms. So, by the usual result of the continuity of integration of a continuous form on a continuous family of cycles (see [B-M 1] ch. IV prop. 2.3.1), we conclude using the   lemma \ref{chgt var.} that the function $\varphi$ is continuous on $T_{0}\setminus \Theta_{0}$.\\
 
 \parag{Step 5} The local boundness on $T$ of the function $\varphi$ is given by the corollary \ref{int 1} which gives the  the estimate $(E)$ by integration.$\hfill \blacksquare$\\
 
 \parag{Remarks} \begin{enumerate}
 \item In the case of a proper family of compact cycles in $X$, it is easy, using results of [B-M 1] chapter IV, to prove that the function $\varphi$ becomes continuous after a suitable modification of the complex space  $T$.
  \item Already in the case of the normalization map, if $\alpha$ is a locally bounded meromorphic function on $X$, the function $x \mapsto \vert \alpha(x)\vert^{2}$ is not continuous on $X$ in general.
  \end{enumerate}

   \section{ Integral dependence  equations and normalized Nash transform}
    We shall use here the notion of linear bundle on a reduced  complex space $X$ . Recall that there is an equivalence between coherent sheaves on $X$ and linear bundles on $X$ which is local and sends a $\mathcal{O}_{X}-$coherent sheaf $\mathcal{F}$  having a local presentation $\mathcal{O}_{X}^{q}\overset{M}{\longrightarrow} \mathcal{O}_{X}^{p} \to \mathcal{F} \to 0$ to the linear bundle $F \subset X\times \C^{p}$ given as the kernel of the map $ ^{t}M : X\times\C^{p} \to X \times \C^{q}$.\\
   The inverse correspondence associates to the linear bundle $F \to X$ the sheaf of linear bundle homomorphism $U \mapsto Hom_{U}(F, \C)$. In particular, when $F$ is a vector bundle, the associated coherent sheaf is the sheaf of holomorphic section of the dual vector bundle $F^{*}$.\\

   In this section we consider a normal complex space $X$.
   
    \subsection{Integral dependence equations}
   
  \parag{Two examples}  \begin{enumerate}
  \item We shall show in section 6.2 that for $k\geq 2$ and $k-1 \geq q \geq k/2$ the form $\omega_{q} := z^{q}.(dx/x - dy/y)$ is a section of the sheaf $\alpha_{S_{k}}^{1}$ where 
  $$S_{k}:= \{(x, y, z) \in \C^{3} \ / \  x.y = z^{k}\}$$ 
  which are not sections of the sheaf $\Omega_{S_{k}}^{1}/torsion$.\\
   But as we have $dx/x + dy/y = k.dz/z$ on $S_{k}$ we obtain the equality
  $$ \omega_{q}^{2} = k.z^{2q - 2}.(dz)^{2} - 4z^{2q-k}.dx.dy ; $$
  so  $\omega_{q}^{2}$ is, for $q \geq k/2$, a section of $S^{2}(\Omega_{S_{k}}^{1})$, the piece of degree $2$ in the symmetric algebra of the sheaf $\Omega_{S_{k}}^{1}$.
  \item We shall show in section 6.4 that on $X := \{(x, y, u, v) \in \C^{4} \ / \ x.y = u.v \}$ the form $a := u.dv\wedge dx/x$ is a section of the sheaf $\alpha_{X}^{2}$ which is not in $\Omega_{X}^{2}/torsion$. But using the following identities on $X$:
  \begin{align*}
  &  u.dv\wedge dx/x + u.dv\wedge dy/y = dv\wedge du \\
  &   u.dv\wedge dy/y + v.du\wedge dy/y = dx\wedge dy
  \end{align*}
  we obtain that 
  $$ a^{2} + a.( du\wedge dv + dx\wedge dy) - (dv\wedge dx).(du\wedge dy) = 0 $$
  which is a homogeneous integral dependence equation for $a$ on the symmetric algebra of the sheaf $\Omega_{X}^{2}/torsion$.
  \end{enumerate}
   
   The following proposition shows that these examples are special cases of a general phenomenon.
   
   \begin{prop}\label{int. dep.}
   Let $X$ be a normal complex space. Then for each integer $q$ the sheaf $\alpha_{X}^{q}$ is the sub-sheaf of meromorphic sections of the sheaf $\Omega_{X}^{q}/torsion$ which satisfy a homogeneous  integral dependence equation over the sheaf $S^{\bullet}(\Omega_{X}^{q})$, the symmetric algebra of the sheaf $ \Omega_{X}^{q}/torsion$.
   \end{prop}
   
   \parag{Proof} Let $\mathcal{T}^{q} \to X$ be the $q-$th Zariski tangent linear bundle over $X$, which is, by definition, the linear bundle associated to the $\mathcal{O}_{X}-$coherent sheaf $\Omega_{X}^{q}$. Then a section $\omega$ of the sheaf $\alpha_{X}^{q}$ defines a morphism of linear bundles $pr : \mathcal{T}^{q} \to X \times \C$ over the open set  $X \setminus S$. As a function on $\mathcal{T}^{q} $ it is locally bounded along $pr^{-1}(S)$ thanks to property $(B)$ in the theorem  \ref{equiv.}. So it defines a meromorphic locally bounded function on the complex space $\mathcal{T}^{q}$ and then satisfies an integral dependence equation of the form
   $$ \omega^{m} + \sum_{j=1}^{m} \ a_{j}.\omega^{m-j} = 0 $$
locally   in an open neighbourhood of $S\times \{0\}$ in $\mathcal{T}^{q}$. As $\omega$ is linear on the fibres, we may assume that the holomorphic function $a_{j}$ is a homogeneous polynomial of degree $j$ along the fibres of $pr : \mathcal{T}^{q} \to X$, and this conclude the proof, because the symmetric algebra of the sheaf $ \Omega_{X}^{q}/torsion$ is the algebra of holomorphic functions on $\mathcal{T}^{q}$ which are polynomial on the fibres modulo these which vanish on $pr^{-1}(X \setminus S)$.$\hfill \blacksquare$\\

   \subsection{Normalized Nash transform}
   
   \parag{Notation} For integers $n < N$ we shall denote $Gr(n, N)$ the grassmannian manifold of sub-vector spaces in $\C^{N}$ of dimension $n$. \\

   Let $X$ be a reduced complex space pure of dimension $n$ and let $S$ its singular locus. Assuming that $X$ is embedded in an open set $U$ in $\C^{N}$ we have a holomorphic map
   $$ \theta : X \setminus S \to Gr(n, N) $$
   sending each point $x \in X \setminus S$ to the $n-$dimensional vector sub-space of $\C^{N}$ which directs the tangent space at $x$ to $X$. This map is holomorphic on $X \setminus S$ and meromorphic along $S$. So the closure of its graph in $X \times Gr(n, N)$ is a proper modification of $X$. We shall note $\mathcal{N} : \hat{X} \to X$ the projection on $X$ of the normalization of this graph. We shall call the (local) {\bf normalized Nash transform} of $X$ this modification. \\
      Let $\pi : \mathcal{U} \to Gr(n, N)$ the universal $n-$vector bundle of $Gr(n, N)$ and let $\mathcal{L}^{q}$ be the sheaf of section of the dual vector bundle to $\Lambda^{q}(\mathcal{U})$. Let $pr : \hat{X} \to Gr(n, N)$ be the projection.
   
   \begin{prop}\label{Nash 1}
   For each integer $q$ there is a canonical isomorphism
    $$ c^{q}: \mathcal{N}^{*}(\alpha_{X}^{q})/torsion \to pr^{*}(\mathcal{L}^{q}) .$$
   \end{prop}
   
   \parag{Proof} Remark first that the pull-back fo holomorphic forms
   $$ \mathcal{N}^{*} : \mathcal{N}^{*}(\Omega_{X}^{\bullet}/torsion) \to \Omega_{\hat{X}}^{\bullet}/torsion $$
   factorizes via the natural sheaf morphisms
   $$ c^{\bullet}: \mathcal{N}^{*}(\Omega_{X}^{\bullet}/torsion) \to \mathcal{L}^{\bullet} \quad {\rm and} \quad \mathcal{L}^{\bullet} \to \Omega_{\hat{X}}^{\bullet}/torsion;$$
   This is a consequence of the fact that $\mathcal{L}^{\bullet}$ is the quotient $ \mathcal{N}^{*}(\Omega_{X}^{\bullet}/torsion)/torsion$.\\
   Consider now a desingularization $\tau : \tilde{X} \to \hat{X}$. Then $\tilde{\tau} := \tau\circ \mathcal{N}$ is a desingularization of $X$ and we have 
   $$ \tilde{\tau}^{**}(\Omega^{\bullet}/torsion) = \tau^{**}(\mathcal{L}^{\bullet}) \subset \Omega_{\tilde{X}}^{\bullet} .$$
   This shows that $\mathcal{N}_{*}(\mathcal{L}^{\bullet})$ is a sub-sheaf of the sheaf $L_{X}^{\bullet} \subset \omega_{X}^{\bullet}$. Moreover, if $\sigma$ is a section of $\mathcal{N}_{*}(\mathcal{L}^{\bullet})$, then $\tilde{\tau}^{*}(\sigma)$ is a section of $\tilde{\tau}^{**}(\Omega_{X}^{\bullet}/torsion)$ thanks to the lemma \ref{pull-back 0} and then $\sigma$ is a section of the sheaf $\alpha_{X}^{\bullet}$.$\hfill \blacksquare$\\
      
   As a consequence of this proposition we obtain that for a normal complex space we have $\alpha_{X}^{q} \simeq \mathcal{N}_{*}(\mathcal{L}^{q})$ for any integer $q \geq 0$. \\
   
   \begin{lemma}\label{increase}
   Let $X$ be a reduced complex space and let $\tau : \tilde{X} \to X$ be any (proper) modification. Then we have a natural inclusion $\alpha_{X}^{\bullet} \hookrightarrow \tau_{*}(\alpha_{\tilde{X}}^{\bullet})$.
   \end{lemma}
   
   \parag{Proof} Consider a desingularization $\theta : \tilde{\tilde{X}} \to \tilde{X}$ and remark that $\pi := \theta\circ \tau$ is a desingularization of $X$. Then we have
   $$ \alpha_{X}^{\bullet} = \pi_{*}(\pi^{**}(\Omega_{X}^{\bullet}) = \tau_{*}(\theta_{*}\big(\theta^{**}(\tau^{**}(\Omega_{X}^{\bullet}))\big) .$$
   Now the equality  $\alpha_{\tilde{X}}^{\bullet} = \theta_{*}(\theta^{**}(\Omega_{\tilde{X}}^{\bullet}))$ and the inclusion  $\tau^{**}(\Omega_{X}^{\bullet}) \subset \Omega_{\tilde{X}}^{\bullet}/torsion $ give
   \begin{align*}
   & \theta^{**}(\tau^{**}(\Omega_{X}^{\bullet}) \subset \theta^{**}(\Omega_{\tilde{X}}^{\bullet}) \\
   & \theta_{*}\big( \theta^{**}(\tau^{**}(\Omega_{X}^{\bullet})\big) \subset  \alpha_{\tilde{X}}^{\bullet} \quad {\rm and \ then} \\
   & \alpha_{X}^{\bullet} \subset \tau_{*}( \alpha_{\tilde{X}}^{\bullet} )
   \end{align*}
   concluding the proof.$\hfill \blacksquare$\\

    \parag{Remark} This shows that when we consider  a sequence of successive modifications over a reduced complex space $X$, the sequence of coherent sub-sheaves $(\tau_{\nu})_{*}(\alpha_{X_{\nu}}^{\bullet})$ is locally stationary on $X$. For instance, this is the case for iterated normalized Nash transforms over a given $X$.\\

\section{Some examples}

 \subsection{Computation of $\omega_{X}^{\bullet}$ for hypersurfaces}
 
 We shall need the following elementary lemma.
 
 \begin{lemma}\label{omega hyp.}
 Let $U$ be an open polydisc in $\C^{n}$ and $D$ an open disc in $\C$. Let $X \subset U \times D$ be a reduced multiform graph of degree $k$ in $U\times D$ with canonical equation $P \in \mathcal{O}(U)[z]$, which is a monic degree $k$ polynomial in $z$. Then we have the inclusion
  $$\Gamma(X, \omega_{X}^{q}) \subset \oplus_{j = 0}^{k-1} \ \frac{z^{j}}{P'(z)}.\Gamma(U, \Omega_{U}^{q})$$
  with equality for $q = n$.
  \end{lemma}
  
  \parag{Proof} First will shall prove the following formula, where $(j, h) \in [0, k]^{2}$:
  $$ det_{j,h}\big[Trace_{X/U}(\frac{z^{j+h}}{P'(z)})\big] = (-1)^{k.(k-1)/2}.$$
  Assume, without loss of generality, that $D$ is centered at the origin with radius $R$. Then for $r > R$  we have, thanks to Cauchy' formula
  $$ Trace_{X/U}(\frac{z^{m}}{P'(z)}) = \frac{1}{2i\pi}.\int_{\vert z\vert = r} \ \frac{z^{m}.dz}{P(z)}. $$
  Then for $m \leq k-2$ put $z = r.e^{i.\theta}$ we obtain
  $$  Trace_{X/U}(\frac{z^{m}}{P'(z)}) = \frac{1}{2\pi}.\int_{0}^{2\pi} \ \frac{r^{m+1-k}.e^{i.(m+1-k)}.d\theta}{1 + O(1/r)} $$
  and letting $r \to +\infty$ gives $0$.
  For $m = k-1$ the same computation gives
  $$  Trace_{X/U}(\frac{z^{k-1}}{P'(z)}) = \frac{1}{2\pi}.\int_{0}^{2\pi} \  \frac{d\theta}{1 + Q((1/r).e^{-i.\theta})} $$
  Where $Q$ is a  polynomial without constant term in $(1/r).e^{-i.\theta}$. So we obtain that $ Trace_{X/U}(\frac{z^{k-1}}{P'(z)}) = 1$. This is enough to obtain the formula above.\\
  
  To prove the inclusion   $\Gamma(X, \omega_{X}^{q}) \subset \frac{1}{P'(z)}.\Gamma(U, \Omega_{U}^{q})$ take $\alpha \in \omega_{X}^{q}$ and write 
  $$\alpha = \sum_{\vert H\vert = q} \ g_{H}.dt^{H}$$
  where $g_{H}$ are degree $\leq k-1$ polynomials in $z$ with meromorphic functions on $U$ as coefficients. As we have $P'(z).dz = - \sum_{h=1}^{n} \ \frac{\partial P}{\partial t_{h}}.dt_{h}$ on $X$, this is possible. Now for any $f \in \mathcal{O}(X)$ we have $Trace_{X/U}[f.\alpha] \in \Omega^{q}(U)$ and this implies that for any $H \subset [1,n], Trace_{X/U}[f.g_{H}] $ is in $\mathcal{O}(U)$. Let $g$ be a meromorphic function on $X$ and assume that we write
  $$ g = \sum_{j=0}^{k-1} \  a_{j}.\frac{z^{j}}{P'(z)} $$
  where $a_{j}, j \in [0, k-1]$ is a meromorphic function on $U$. This is always possible for the $g_{H}$  as we can see in what follows. Let $m_{p}:= Trace_{X/U}[z^{p}.g]$ for $p \in [0, k-1]$. Then we have the linear system in the $(a_{j}), j \in [0, k-1]$:
  $$ \sum_{j=0}^{k-1} \ a_{j}.Trace_{X/U}[\frac{z^{p+j}}{P'(z)}] = m_{p} \quad \forall p \in [0, k-1].$$
  But the determinant of this linear system is $ (-1)^{k.(k-1)/2}$, so this implies, if we assume that the functions $m_{p}$ are holomorphic on $U$, that the functions  $a_{j}$ for $ j \in [0, k-1]$,  are holomorphic in $U$ and so that $g$ is in $\frac{1}{P'(z)}.\mathcal{O}(X)$. Then our inclusion is proved, as $\mathcal{O}(X) = \sum_{j=0}^{k-1} \ \mathcal{O}(U).z^{j}$.\\
  Note that in the situation above, the condition in order that $\alpha = \frac{1}{P'(z)}.\Omega^{q}(U)$ will be in $\omega^{q}(X)$ is that for any $j \in [0, k-1]$ the $(q+1)-$forms
  $$ Trace_{X/U}[z^{j}.dz\wedge \alpha]$$
   are holomorphic for all $j \in [0, k-1].$
  This is consequence of the fact that  for any $\beta \in \Omega^{p}(X)$  the $(p+1)-$form $Trace_{X/U}[\alpha\wedge \beta]$ must be holomorphic (see [B. 78] for this characterization of the sheaf $\omega_{X}^{\bullet}$).  For $q = n$ this extra condition is empty, so the equality occurs.  $\hfill \blacksquare$\\
  
  \parag{Remark} For a general reduced multiform graph $X \subset U \times B$ where $B$ is now a polydisc in $\C^{p}$, for any linear form $l$ in $\C^{p}$ which separates generically the fibers of the projection $\pi : X \to U$, the map $id_{U}\times l : U\times B \to U\times \C$ is proper and generically injective on $X$. If we define $Y_{l} := (id_{U}\times l)(X)$, we are in the situation of the lemma above, and, as the direct image by $\pi$ induces an injective sheaf map $\pi_{*} : \omega_{X}^{\bullet} \to \pi^{*}(\omega_{Y_{l}}^{\bullet})$, we obtain the inclusion
   $$\omega_{X}^{\bullet} \subset \oplus_{j=0}^{k-1} \  \frac{l(x)^{j}}{P'_{l}}.\Omega_{U}^{\bullet}$$
   for any such $l$, where $P_{l}$ is the canonical equation for $Y_{l}$ (see [B-M 1] chapter II). Note that the canonical equation $P_{l}$ is obtained from the canonical equation of  the reduced multiform graph $X$ by the evaluation at $l$ (with $z = l(x)$); see {\it loc. cit.}.$\hfill \square$\\
  
  Note that, if $X$ is a reduced complex space of pure  dimension $n$,  a section $\alpha \in \omega_{X}^{n}$ is in $L_{X}^{n}$ iff $\alpha\wedge \bar \alpha$ is locally integrable on $X$. The analogous characterization for $p < n$, involves local integrability of $\alpha\wedge \bar\alpha$ on  all $p-$dimensional irreducible analytic subset $Y \subset X$ not contained in the singular set of $X$; so it may be useful as a necessary condition but very difficult to check as a sufficient condition.

 \subsection{The case $X := \{(x,y,z) \in \C^{3}\ / \ x.y = z^{k}\}, k \geq 2$}
 
 \parag{Notation} After blow-up $(x, y, z)$ in $\C^{3}$ the homogeneous coordinates in $\mathbb{P}_{2}$ will be $(\alpha, \beta, \gamma)$. The symetry between $x$ and $y$ allows to consider only the chart $\{\alpha \not= 0\}$ on which we put $b := \beta/\alpha, c := \gamma/\alpha$ and the chart $\{\gamma \not= 0\}$ on which we put $a := \alpha/\gamma, b := \beta/\gamma.$\\
 
 Our first example will be the normal complex spaces, where $k \in \mathbb{N}, k \geq 2$
  $$ X := S_{k} := \{(x, y, z) \in \C^{3} \ / \  x.y = z^{k} \}.$$
  Note that $S_{0}$ and $S_{1}$ are smooth complex surfaces.\\
  
  \begin{lemma}\label{nearly smooth}
  For any $k \geq 2$ the normal complex space $S_{k}$ is nearly smooth\footnote{See [B-M. 17]}. So we have $L_{S_{k}}^{\bullet} = \omega_{S_{k}}^{\bullet}$ for any $k$.
  \end{lemma}
  
  \parag{Proof} Let $\zeta$ be a $k-$th primitive root of $1$. Then $S_{k}$ is somorphic to the quotient of $\C^{2}$ by the action of the automorphism  $\theta(u,v) = (\zeta.u, \zeta^{-1}.v)$. The quotient map is given by $q(u, v) = (u^{k}, v^{k}, u.v) \in \C^{3}$.$\hfill \blacksquare$\\ 
  
Now compute the sheaf $\omega_{X}^{h}$ for $h \in [0, 2]$. We have $\omega_{X}^{0} = \mathcal{O}_{X}$ as $X$ is normal, and $\omega_{X}^{2} = \mathcal{O}_{X}.\frac{dx\wedge dy}{z^{k-1}}$. A rather easy computation shows that the quotient $\omega_{X}^{1}\big/\Omega_{X}^{1}$ is generated on $\mathcal{O}_{X}$  by the image of $x.dy\big/z^{k-1} = -y.dx/z^{k-1} + k.dz$ which is annihilated in this quotient  by $x, y$ and $z^{k-1}$.\\
 
 \begin{lemma}\label{L2+}
 For any $k \geq 2$  the sheaf  $\alpha_{S_{k}}^{2}$ coincides with $\Omega_{S_{k}}^{2}\big/torsion$.
 \end{lemma}
 
 \parag{Proof}  Remark that for $k = 0, 1$ the lemma is obvious as $S_{k}$ is smooth. We shall prove the lemma by induction on $k \geq 2$. \\ We have to consider the case $k=2$ first because it appears that the computation is special in this case (see the denominator $k-2$ in the computation for $k \geq 3$).\\
 
 For $k = 2$ after blowing-up the origin we have a smooth manifold:\\
  In the chart $\{\alpha \not= 0\}$ we have $y = x.b,\quad z = x.c, \quad b = c^{2} \quad$  so $(x,c)$ is a coordinate system in this chart and 
 $$ \frac{dx \wedge dy}{z} = 2.dx\wedge dc $$
 is holomorphic but not in $\tau^{*}(\Omega_{S_{2}}^{2}) \simeq \mathcal{O}_{\tilde{X}}.x.dx\wedge dc$.\\
 In the chart $\{\gamma \not= 0\}$ we have $x = z.a,\quad  y = z.b,\quad  a.b = 1$. So $(z, a)$ is a coordinate system with $a \not= 0$. Then
  $$ \frac{dx \wedge dy}{z} = -2dz\wedge da/a $$
  which is holomorphic but not in $\tau^{**}(\Omega_{S_{2}}^{2}) \simeq \mathcal{O}_{\tilde{X}}.z.dz\wedge da$.\\
  The assertion is proved for $k=2$.\\

    As the assertion is proved for  $k =  2,$ we may assume that, for $k \geq 3$  the equality is proved for $S_{k-2}$. Then let $\tilde{X} \to X := S_{k}$ be the blow-up of $S_{k}$ at the singular point $x = y = z = 0$.  In the chart $\{\gamma \not= 0\}$ of $\tilde{X}$ we have the relations
 $$ x = a.z, \quad y = b.z \quad a.b = z^{k-2} \quad $$
 and we find a copy of $S_{k-2}$.  For $k \geq 3$ we have
 \begin{align*}
 & dx\wedge dy = \frac{k}{k-2}.z^{2}.da\wedge db = k.\frac{a.b}{k-2}.\frac{da\wedge db}{z^{k-4}} ,\\ 
 & \quad  dx\wedge dz = \frac{a}{k-2}.\frac{da\wedge db}{z^{k-4}}, \quad  dy\wedge dz = \frac{b}{k-2}.\frac{da\wedge db}{z^{k-4}}.
 \end{align*} 
 So in this chart
  $$\tau^{**}(\Omega_{S_{k}}^{2}\big/torsion) = \mathcal{O}_{S_{k-2}}.\big(a, b\big).\frac{da\wedge db}{z^{k-4}}$$
 and, as a  consequence of the fact that $z^{k-q-2}$ is not in the ideal $\big(a, b\big).\mathcal{O}_{S_{k-2}}$ for $q \geq 1$, for no $q \geq 1$ the $2-$form $dx\wedge dy\big/z^{q}$ is a section of the sheaf $\tau^{**}(\Omega_{S_{k}}^{2}\big/torsion) $ near the origin $a = b = z = 0$ in this chart.  So the sheaf $\alpha_{S_{k}}^{2}$ is equal to $\Omega_{S_{k}}^{2}\big/torsion$.$\hfill \blacksquare$\\
          
  \begin{lemma}\label{k geq 0}
 For all $k \geq 0$  the vector space $L_{S_{k}}^{1}\big/\alpha_{S_{k}}^{1}$ has  dimension $p = [(k-1)/2] $  the integral part  of  $(k-1)/2$. A basis is given by the $1-$forms $x.dy/z^{q}$ for $q$ in  $[ [k/2]+1, k-1]$, for $k \geq 3$.
 \end{lemma}

 \parag{Proof} The lemma is clear for $k = 0, 1$.  We have seen that after blowing-up the singular point in $S_{k}$ for any $k \geq 2$ we find only one singular point of the type $S_{k-2}$ in the chart $\{\gamma \not= 0\}$ and that the form $x.dy/z^{q}$ is given by the following computation in this chart $\{\gamma \not= 0\}$ :
 $$ x = z.a, \quad y = z.b,\quad  a.b = z^{k-2} \quad x.dy/z^{q} = a.db/z^{q-2} + z^{k-q-1}.dz .$$
 So assuming that $ k = 2p+1 \geq 3$ and that we know that the forms $a.db/z^{q-2}$ for $q-2$ in $\{ p, \dots, 2p-2 \}$ is a basis of the vector space $L_{S_{2p-1}}^{1}\big/\alpha_{S_{2p-1}}^{1}$, this  implies  that the forms $x.dy/z^{q}$ for $q = p+2, \dots, 2p$ is a basis of a vector space of dimension $p-1$ in $L_{S_{2p+1}}^{1}\big/\alpha_{S_{2p+1}}^{1}$ because we have a linear  map $L_{S_{2p+1}}^{1}/\alpha_{S_{2p-1}}^{1}/\alpha_{S_{2p-1}}^{1}$ given by the pull-back. It is then enough to prove that for $q = p+1$ the form $x.dy/z^{p+1}$ is not in $\alpha_{S_{2p+1}}^{1}$ and that it  completes in a basis the previous free system in $L_{S_{2p+1}}^{1}\big/\alpha_{S_{2p+1}}^{1}$ to conclude. In the last chart $\{\gamma \not= 0\}$ in the desingularization process of $S_{2p+1}$ by blowing up the unique singular point at each step, we reach the following relations:
 $$ x= u^{p}.v^{p+1}, y = u^{p+1}.v^{p}, z = u.v \quad x.dy\big/z^{p+1} = (p+1).u^{p-1}.v^{p}.du + p.u^{p}.v^{p-1}.dv $$
 where $(u, v) \in \C^{2}$ is a local coordinate system.\\
 Now assume that there exists holomorphic functions $\lambda, \mu, \nu$ near the origin on $\C^{2}$ and a polynomial $Q(z) := 1+  \sum_{j=1}^{p} q_{j}.z^{j}$  such that we have
  $$Q(z).x.dy/z^{p+1} = \lambda(u,v).dx + \mu(u,v).dy + \nu(u,v).dz$$
  near $u = v = 0$. This implies the equalities
 \begin{align*}
 & Q(u.v).(p+1).u^{p-1}.v^{p-1} = \lambda.p.u^{p-1}.v^{p} + \mu.(p+1).u^{p}.v^{p-1} + \nu \\
 &  Q(u.v).p.u^{p-1}.v^{p-1} =  \lambda.(p+1).u^{p-1}.v^{p} +  \mu.p.u^{p}.v^{p-1} + \nu \quad {\rm and \ then } \\
 &  Q(u.v).u^{p-1}.v^{p-1} = -  \lambda.u^{p-1}.v^{p} + \mu.u^{p}.v^{p-1}  \quad  {\rm and \ so} \\
 & Q(u.v) = -\lambda.v + \mu.u
 \end{align*}
 which is a contradiction. Note that the same computation gives
 \begin{align*}
 &  x.dy/z^{p} = (p+1).u^{p}.v^{p+1}.du + p.u^{p+1}.v^{p}.dv \\
 & \qquad = (p+1).x.du + p.y.dv \quad  {\rm and \ as} \  xu = yv = z^{p+1} \\
 & \qquad = (p+1).(dz^{p+1} - u.dx) + p.(dz^{p+1} - v.dy) \\
 & \qquad= (p+1).dz^{p+1} - (p+1).u.dx - p.v.dy \in \tau^{**}(\Omega^{1}_{S_{2p+1}}\big/torsion) \\
 \end{align*}
 
Now assume that $k = 2p$ with $p \geq 2$ then in the last chart $\{\gamma \not= 0\}$ we shall have, with coordinates $(z, u)$ with $u \not= 0$
$$  x= z^{p}.u,\quad y = z^{p}/u  \quad {\rm so} \quad   x.dy/z^{p+1} = p.z^{p-2}.dz - z^{p-1}.du/u .$$
Then assume that, again with $Q(z) := 1 +  \sum_{j=1}^{p} q_{j}.z^{j}$
$$Q(z).z^{p-1}.du = \lambda.(z^{p}.du + p.u.z^{p-1}.dz) +  \mu.(-z^{p}.du/u^{2} + p.u^{-1}.z^{p-1}.dz) + \nu.dz ; $$ 
 it implies that 
$$ Q(z).z^{p-1} = \lambda.z^{p} - \mu.z^{p}.u^{-2} \quad {\rm and } \quad  Q(z) = (\lambda - \mu/u^{2}).z$$
which gives a contradiction for $z = 0$.\\
 Note that  it easy to see that $z^{p}.du = dx - p.u.z^{p-1}.dz \in \tau^{**}(\Omega_{S_{2p}}^{1}\big/torsion).\hfill \blacksquare$\\
 
 \newpage
  
  \subsection{ The case $X := \{ (x,y,z) \in \C^{3}\ / \  x^{3} + y^{3} + z^{3} = 0 \}$}
  
  Now consider $X := \{ (x, y, z) \in \C^{3} \ / \  x^{3} + y^{3} + z^{3} = 0 \}$. The lemma \ref{omega hyp.} gives the inclusion
   $$\omega_{X}^{1}\subset \frac{1}{z^{2}}.\Omega_{\C^{2}}^{1} $$
   where $x, y$ are the coordinates on $\C^{2}$. An easy computation shows that the form $\alpha := (x.dy - y.dx)/z^{2}$ generates $\omega_{X}^{1}\big/\Omega_{X}^{1}$.\\
   
     Let $\tau : \tilde{X} \to X$ the blowing-up at the origin of $X$ is a desingularization. In the chart. $\{\gamma \not= 0\}$ let $a := \alpha/\gamma$ and $b := \beta/\gamma$; then  we have the relations 
      $$x = z.a,\quad  y = z.b,\quad  a^{3} + b^{3} + 1 = 0$$
      and then we can choose $(z, a)$ or $(z, b)$  as local  coordinates. Then we have
     $$ \alpha = a.db - b.da = db/a^{2} = -da/b^{2}.$$
     In the chart $\{\alpha \not= 0\}$ we have  $$y = u.x, z = v.x, u^{3} + v^{3} + 1 = 0.$$
      Then we can choose $(x, u)$ or $(x, v)$ as local coordinates and  $\alpha = du/v^{2} = -dv/u^{2}$. 
     This shows that $ \omega_{X}^{1} = L_{X}^{1} $. But $\alpha$ does not vanish on the exceptional divisor, so $\alpha$ is not a section of $\alpha_{X}^{1}$. \\
     But,   in the first chart, 
     $$z.\alpha = z.a.db - z.b.da = a.dy - a.b.dz - b.dx + a.b.dz = a.dy - b.dx \in \tau^{**}(\Omega_{X}^{1})$$
    and  in the second chart 
     $$x.\alpha = x.du/v^{2} = dy/v^{2}-u.dx/v^{2}= -dz/u^{2} + v.dx/u^{2}$$
      also belong to $ \tau^{**}(\Omega_{X}^{1})$.\\
     Then $x.\alpha, y.\alpha$ and $z.\alpha$ are  sections of $\alpha_{X}^{1}$ and  the quotient $L_{X}^{1}\big/\alpha_{X}^{1}$ is a vector space of dimension $1$ with basis $\alpha $.$\hfill \blacksquare$\\
     
     Note that $x.y.z.\alpha$ is not a section of $\Omega_{X}^{1}/torsion$ because if this is not the case, we can write
     $$  x.y.(x.dy- y.dx)= z.\big[\lambda.dx  + \mu. dy + \nu.dz + \rho.df + \sigma.f\big] $$
     in $\C^{3}$, where $\lambda, \mu, \nu$ where homogeneous of degree $2$,  $\rho$ is a complex number and where $\sigma := u.dx + v.dy + w.dz$ with $u, v, w$ 
      complex numbers. This gives, for instance $-x.y^{2} = z.\lambda + 3z.\rho.x^{2} + u.f$ which is impossible.\\
      So the vector space $\alpha_{X}^{1}/\Omega_{X}^{1}$ has dimension at least $2$. The complete determination of the quotient $\alpha_{X}^{1}/\Omega_{X}^{1}$  is a non trivial exercise left to the reader...

      \bigskip
     
     \begin{lemma}\label{tout distinct}
     For $X := \{ (x,y,z) \in \C^{3}\ / \  x^{3} + y^{3} + z^{3} = 0 \}$  we have
     $$dim_{\C}\  \alpha_{X}^{2}\big/\Omega_{X}^{2} = 2, \quad dim_{\C} \ L_{X}^{2}\big/\alpha_{X}^{2} = 3 \quad dim_{\C}\ \omega_{X}^{2}\big/L_{X}^{2} = 1 .$$
     \end{lemma}
     
     \parag{Proof} After blowing-up $(x, y, z)$ in $\C^{3}$ we consider the chart $\{\gamma \not= 0\}$ as above. We have
     $$ \omega := \frac{dx\wedge dy}{z^{2}} = -\frac{dz}{z}\wedge\frac{db}{a^{2}} =  \frac{dz}{z}\wedge\frac{da}{b^{2}} .$$
     Then $x.\omega, y.\omega, z\omega$ are holomorphic in this chart, as we have $x = z.a$ and $y = z.b$ and  this chart is enough as $dx\wedge dy/z^{2} = dy\wedge dz\big/x^{2} = dz\wedge dx\big/y^{2}$
     so $x.\omega, y.\omega, z.\omega$ belongs to $L_{X}^{2}$. \\
     But this is not the case for $\omega$. So $\dim \omega^{2}/L_{X}^{2} = 1$.\\
     The sheaf $\tau^{**}(\Omega_{X}^{2}\big/torsion)$ in this chart is generted by
        $$z.(da/a^{2})\wedge dz = -z.(db/b^{2})\wedge dz.$$
        Then it is equal to $z.\Omega_{\tilde{X}}^{2}$ in this chart. So a section in $L_{X}^{2}$ is in $\alpha_{X}^{2}$ if and only if it belongs to $(x.L_{X}^{2})\cap (y.L_{X}^{2})\cap(z.L_{X}^{2})$. This intersection is generated by $x.y.\omega, y.z.\omega, z.x.\omega$ as a $\mathcal{O}_{X}-$module. The vector space $L_{X}^{2}/\alpha_{X}^{2}$ is generated by  $x.\omega, y.\omega, z.\omega$ because $x^{2}.\omega, y^{2}.\omega, z^{2}.\omega$ are in $\Omega_{X}^{2}\subset \alpha^{2}_{X}$. We let  to the reader the proof that they give a basis of $L_{X}^{2}/\alpha_{X}^{2}$.\\
                
        Let us prove that $x.y.z.\omega$ is not in $\Omega_{X}^{2}/torsion$.\\
        Assume that  $x.y.z.\omega \in \Omega_{X}^{2}/torsion$.  Then we can write on $\C^{3}$:
        $$ x.y.dx\wedge dy -  z\big[\lambda.dx\wedge dy + \mu.dy\wedge dz + \nu.dz\wedge dx + (a.dx + b.dy + c. dz) \wedge df\big]  = 0$$
        where we can assume that $\lambda, \mu, \nu$ are linear forms on $\C^{3}$ and $a, b, c$ are complex number, using the homogeneity of the situation. The coefficient of $dx\wedge dy$ in this identity is equal to
        $x.y - z.\lambda - a.y^{2} + b.x^{2}$ which cannot be identically zero. Contradiction.\\
        As it is easy to see that $x.y.\omega= y.z.\omega= z.x.\omega$  and $x.y.z.\omega$ are linearly independent over $\C$ (different homogeneities) we conclude that $\dim \alpha_{X}^{2}/\Omega_{X}^{2} = 2$.$\hfill \blacksquare$\\
        
        \parag{Remark} We have on $X$ 
   $$      \omega := \frac{dx\wedge dy}{z^{2}} = \frac{dy\wedge dz}{x^{2}} = \frac{dz\wedge dx}{y^{2}} $$
   so
    \begin{align*}
    &(x.y.\omega)^{2} = \frac{x^{2}.y^{2}.(dx\wedge dy)^{2}}{z^{4}} = \frac{x^{2}.dx\wedge dy}{z^{2}}.\frac{y^{2}.(dx\wedge dy)}{z^{2}} = (dz\wedge dy).(dx \wedge dz),
    \end{align*}
   because on $X$ we have $x^{2}.dx\wedge dy = -z^{2}.dz\wedge dy$ \ and \ $y^{2}.dx\wedge dy = -z^{2}.dx\wedge dz$. This gives an integral dependence relation for $x.y.\omega$ on $\Omega_{X}^{2}/torsion$.\\

     \subsection{The case $X := \{(x, y, u, v) \in \C^{4}\ / \  x.y = u.v \}$ }
 
 \begin{lemma}\label{L3}
 The sheaf $L_{X}^{3}$ is equal to $\omega_{X}^{3}$ and is given by $\mathcal{O}_{X}.\omega$ where we define 
  $$\omega := \frac{dy\wedge du\wedge dv}{y}.$$ 
 Moreover, $\omega$ does not belong to $\alpha_{X}^{3}$.
 \end{lemma}
 
 \parag{proof} On $X$ we have $x.dy + y.dx= u.dv + v.du$   
 $$ \omega = - \frac{dx\wedge du\wedge dv}{x} =  \frac{du\wedge dx \wedge dy}{u} =  \frac{dv \wedge dx \wedge dy}{v} .$$
 To see that $\omega_{X}^{3} = \mathcal{O}_{X}.\omega$ it is enough ($X$ is a hypersurface !) to see  that 
 $$\omega\wedge df/f = dx\wedge dy\wedge du\wedge dv/f$$
  where $f := x.y - u.v$. This is clear.\\
 
 Using the symetries between the coordinates, it is enough to see that $\tau^{*}(\omega)$ is holomorphic in the first chart of the strict transform $\tilde{X}$  of $X$ by the blow-up at the origin in $\C^{4}$ to show that $\omega$ is a section of $L_{X}^{3}$. Let $y = \lambda.x, u = \mu.x, v = \nu.x$. Then 
 $$\tau^{*}(\omega) = -\frac{dx}{x}\wedge x.d\mu \wedge x.d\nu = -x.dx\wedge d\mu \wedge d\nu $$
 where $x, \mu, \nu$ are the coordinates for $\tilde{X}$ in this chart (and $\lambda = \mu.\nu$). So $\omega \in L_{X}^{3}$. \\
 To see that $\omega$ is not in $\alpha_{X}^{3}$ it is enough to see that $\omega$ does not belongs to $\tau^{**}(\Omega_{X}^{3})$ in the first chart above. An easy computation show that $\tau^{**}(\Omega_{X}^{3})$ is generated by $\pi^{**}(dx\wedge du\wedge dv) = x^{2}.dx\wedge d\mu \wedge d\nu $ and so $\omega = -x.dx\wedge d\mu \wedge d\nu$ does not belong to $\tau^{**}(\Omega_{X}^{3})$.$\hfill \blacksquare$\\
 
 \begin{lemma}\label{L2}
 The meromorphic form $w := u.dv\wedge dx/x$ is a section of  $\alpha_{X}^{2}$ but it is not a section of  $\Omega_{X}^{2}/torsion$ and its differential is not a section of  $\alpha_{X}^{3}$.
 \end{lemma}
 
 \parag{proof} As
  $$u.dv\wedge dx/x + v.du\wedge dx/x = dx\wedge dy$$
   is holomorphic on $X$, $u$ and $v$ play the same role for this form modulo holomorphic forms. Also $u.dv \wedge (dx/x + dy/y) = dv\wedge du$ so $x$ and $y$ play also the same role modulo holomorphic forms on $X$. So it is enough to see that in the first chart of the strict transform $\tilde{X}$  of $X$ by the blow-up at the origin in $\C^{4}$ the form $\tau^{**}(w)$ is a section of $\tau^{**}(\Omega_{X}^{2})$ to prove that $w$ is a section of $\alpha_{X}^{2}$. Using the same coordinates as above we obtain
 $$ \tau^{**}(w) = \mu.x.d(\nu.x)\wedge dx/x = \mu.x.d\nu\wedge dx = \mu.dv\wedge dx $$
 which is a section of $\tau^{**}(\Omega_{X}^{2})$.\\
 To see that $w$ is not a section of $\Omega_{X}^{2}/torsion$ assume the contrary. Then, by symetry\footnote{or using $(u.dv + v. du)\wedge dx/x = dy\wedge dx$.} $w' := v.du\wedge dx/x$ is also a section of $\Omega_{X}^{2}/torsion$ and the differential of  $w- w'$ must be a section of  $\Omega_{X}^{3}/torsion$. But we have already seen that $2.\omega = -d(w- w')$ is not a section of $\alpha_{X}^{3}$. Contradiction. $\hfill \blacksquare$\\
 
 Note that an integral dependence relation on the symmetric algebra of the sheaf $\Omega_{X}^{2}/torsion$ for $w$ is given in the second example of the begining of the section 5.1.
 
 \begin{lemma}\label{L1}
 We have $\Omega_{X}^{1}/torsion = \alpha_{X}^{1} = L_{X}^{1} = \omega_{X}^{1}$.
 \end{lemma}
 
 \parag{Proof} Write $X := \{x_{1}^{2} + x_{2}^{2} + x_{3}^{2} + x_{4}^{2} = 0\} \subset \C^{4}$. Then thanks to the lemma \ref{omega hyp.} we have $\omega_{X}^{1} \subset \oplus_{j=0}^{1} \ \frac{x_{4}^{j}}{x_{4}}.\Omega_{\C^{3}}^{1} $. To prove that $\Omega_{X}^{1}/torsion = \omega_{X}^{1}$ it is enough to consider a section  $v := (a.dx_{1}+ b.dx_{2}+ c.dx_{3})/x_{4}$ in $\omega_{X}^{1}$ and to show that it is a section of $\Omega_{X}^{1}/torsion$. But then $Trace_{\pi}(v\wedge dx_{4})$ must be a holomorphic form on $\C^{3}$, where $\pi : X \to \C^{3}$ is the projection which makes $X$ a branched covering of degree $2$. This condition implies $(x_{1}.dx_{1} + x_{2}.dx_{2}+ x_{3}.dx_{3}) \wedge (a.dx_{1}+ b.dx_{2}+ c.dx_{3}) = 0$ and then   $(a.dx_{1}+ b.dx_{2}+ c.dx_{3}) = u.x_{4}.dx_{4}$ on $X$, where $u$ is holomorphic on $\C^{3}$. So $v$ is a section of $\Omega_{X}^{1}/torsion$.$\hfill \blacksquare$\\
 
 \begin{lemma}\label{L2}
 We have $\omega_{X}^{2} = \Omega_{X}^{2}/torsion \oplus \C.\eta $ where
 $$ \eta := \frac{x_{1}.dx_{2}\wedge dx_{3} + x_{2}.dx_{3}\wedge dx_{1} + x_{3}.dx_{1}\wedge dx_{2}}{x_{4}}.$$
 \end{lemma}
 
 \parag{Proof} Write $\omega := (a.dx_{1}\wedge dx_{2} + b.dx_{2}\wedge dx_{3} + c.dx_{3}\wedge dx_{1})/x_{4}$ where $a,b,c$ are holomorphic on $\C^{3}$. Then $\omega$  is in $\omega_{X}^{2}$ if and only if $Trace_{\pi}(dx_{4}\wedge \omega)$ is a section of $\Omega_{\C^{3}}^{3}$. This is satisfyed if and only if $a.x_{3} + b.x_{1} + c.x_{2}$ is a multiple of $\xi := x_{1}^{2}+ x_{2}^{2}+ x_{3}^{2} $ in $\mathcal{O}_{\C^{3}}$. This gives the relation $ (a - g.x_{3}).x_{3} + (b - g.x_{1}).x_{1} + (c - g.x_{2}).x_{2} = 0$. And, as $x_{1}, x_{2}, x_{3}$ is a regular sequence, this implies
 $$ a =  g.x_{3} + \lambda.x_{1} + \mu.x_{2}, \quad b = g.x_{1} + \lambda'.x_{2} - \lambda.x_{3}, \quad c = g.x_{2} - \lambda'.x_{1} - \mu.x_{3} $$
 where $\lambda, \lambda', \mu$ are in  $\mathcal{O}_{\C^{3}}$. This shows that $\omega_{X}^{2}$ is generated as a $\mathcal{O}_{X}-$module  by $\Omega_{X}^{2}$ and $\eta$.
 Note that we already know that $\eta$ is not a section of  $\Omega_{X}^{2}/torsion$ as we have shown that $\omega_{X}^{2}$ is not equal to $\Omega_{X}^{2}/torsion$ and that $x_{i}.\eta$ is in $\Omega_{X}^{2}/torsion$ for $i = 1, 2, 3, 4$ :\\
 for instance
 \begin{align*}
 & \frac{x_{1}.\eta}{x_{4}} = = \frac{x_{1}}{x_{4}}.(x_{1}.dx_{2}\wedge dx_{3} + x_{2}.dx_{3}\wedge dx_{1} + x_{3}.dx_{1}\wedge dx_{2}) \\
 & \quad = \frac{1}{x_{4}}.(-(x_{2}^{2}+ x_{3}^{2} + x_{4}^{2}).dx_{2}\wedge dx_{3} + (x_{2}.dx_{3} - x_{3}.dx_{2})\wedge x_{1}.dx_{1}) \\
 & \quad = \frac{1}{x_{4}}.(-(x_{2}^{2}+ x_{3}^{2} + x_{4}^{2}).dx_{2}\wedge dx_{3} + (x_{2}.dx_{3} - x_{3}.dx_{2})\wedge (-x_{2}.dx_{2}-x_{3}.dx_{3}- x_{4}.dx_{4})) \\
 & \quad =  - x_{4}.dx_{2}\wedge dx_{3}- x_{2}.dx_{3}\wedge dx_{4} + x_{3}.dx_{2}\wedge dx_{4}\\
 \end{align*}
 proving our claim. $\hfill \blacksquare$\\

 \subsection{The case $X := \{(x, y, z, t) \in \C^{4} \ / \  x.y.z = t^{3} \}$}
 
 Remark first that the form $\omega_{1}:= y.z.dx/t^{2}$ is in $\omega_{X}^{1}$ because we have, with the notation $f := x.y.z - t^{3}$:
 $$ \omega_{1}\wedge df= \frac{z.t^{3}.dx\wedge dy + y.t^{3}.dx\wedge dz + 3t^{2}.y.z.dx\wedge dt}{t^{2}} \in \Omega_{\C^{4}}^{2}$$
 modulo $(f/t^{2}).\Omega_{\C^{4}}^{2}$ which allows to conclude as $t$ is not a zero divisor in $X$.\\
 Consider now the following sections of $\omega_{X}^{1}$:
 $$ u := t.\omega_{1} \quad v := t.\omega_{2} \quad w := t.\omega_{3} $$
 where $\omega_{2}$ and $\omega_{3}$ are deduced from $\omega_{1}$ respectively by $x \to y, y \to z, z \to x$ and $x \to z, y \to x, z \to y$. Then we have in the symmetric algebra of $\Omega_{X}^{1}$
 $$ u + v + w = 3t.dt \quad u.v + v.w + w.u = t.(z.dx.dy + x.dy.dz + z.dx.dy) \quad u.v.w = t^{3}.dx.dy.dz .$$
 This shows that $u, v, w$ satisfy the following integral dependence relation over the symmetric algebra of $\Omega_{X}^{1}$
 $$  \Theta^{3} - 3t.dt.\Theta^{2} + t.(z.dx.dy + x.dy.dz + z.dx.dy).\Theta -  t^{3}.dx.dy.dz = 0 .$$
 
 Note that, because the coefficient of $\Theta$ does not belong to $(t^{2})$, we do not obtain an integral dependence relation over the symmetric algebra of $\Omega_{X}^{1}$  for $\Theta/t$ so for the forms $\omega_{i}, i = 1, 2, 3$ ! In fact they are not sections of the sheaf $\alpha_{X}^{1}$ (for instance the restriction of $\omega_{1}$ to the surface $S_{3} \simeq \{z = 1\} \cap X$ is not in $\alpha_{S_{3}}^{1}$ (see sub-section 6.2).\\
 
  Let us now verify that $t.u$ is not a section of $\Omega_{X}^{1}/torsion$. Let assume that we can write
 $$ y.z.dx = t.\big(\lambda.dx + \mu.dy + \nu.dz + \theta.dt) \quad {\rm modulo} \ f.\Omega_{X}^{1} + \mathcal{O}_{X}.df $$
 then, by homogeneity, we may assume that $\lambda, \mu, \nu$ are homogeneous of degree $2$ and 
 $$ y.z.dx = t.\big(\lambda.dx + \mu.dy + \nu.dz + \theta.dt) + \sigma.df$$
 where $\sigma$ is a constant. This implies
 $$ y.z.(1 - \sigma) - t.\lambda  = 0, \quad  t.\mu + \sigma.x.z = 0   $$
 which is already enough to obtain a contradiction, as these equations imply $\sigma = 1$ and $\sigma = 0$ respectively.$\hfill \blacksquare$\\

 \parag{Remark} Using the map $((x, y, z) \mapsto (x+y, x+j.y, x+j^{2}.y, -z)$  which sends the previous $Y := \{x^{3}+ y^{3} + z^{3} = 0\}$ to $X = \{  x.y.z = t^{3} \}$ allows to find an integral equation over the symmetric algebra of $\Omega_{Y}^{1}$ of the section
  $$\frac{(x^{2}+ y^{2} -x.y).d(x+y)}{z}$$
  of $\alpha_{Y}^{1}$.

 \bigskip
 
 \parag{References}
 
 \begin{itemize}
 
 \item{[B.78]} Barlet, D. {\it Le faisceau $\omega_{X}^{\bullet}$ on a reduced complex space},  in Sem. F. Norguet III, Lecture Notes, vol. 670, Springer Verlag (1978), p. 187-204.
 
 \item{[B-M 1]} Barlet, D. et Magnusson, J. {\it Cycles analytiques complexes I : th\'eor\`{e}mes de pr\'eparation des cycles}, Cours Sp\'ecialis\'es 22, Soci\'et\'e Math\'ematiques de France, Paris 2014.
 \item{[B.15]}  Barlet, D. {\it Strongly quasi-proper maps and the f-flattning theorem},  arXiv:1504.01579 (math.CV)
 
 \item{[B-M.17]} Barlet, D. et Magnusson, J. {\it Nearly-smooth complex spaces} to appear soon in math-arXiv.

 \end{itemize}

\end{document}